\documentclass[reqno,12pt]{amsart}
\usepackage{ot1patch}
\usepackage{amsmath,amsfonts,amsthm,amssymb,enumerate,graphicx,esint}
\usepackage{caption}
\usepackage[dvipsnames]{xcolor}
\usepackage[margin=1in]{geometry}
\usepackage{tikz,tikz-cd,mathtools,soul}

\setstcolor{red}

\let\div\relax
\DeclareMathOperator\div{div}

\DeclareMathOperator\curl{curl}

\DeclareMathOperator\supp{supp}
\DeclareMathOperator\tr{tr}
\DeclareMathOperator\dist{dist}
\DeclareMathOperator*\osc{osc}
\newcommand\dd[0]{\partial}
\newcommand\grad[0]{\nabla}
\newcommand\ulin[0]{u^\flat}
\newcommand\unlin[0]{u^\sharp}

\newcommand\eq[1]{\begin{align}{#1}\end{align}}
\newcommand\eqn[1]{\begin{align*}{#1}\end{align*}}
\newcommand\Rd[0]{{\mathbb R^3}}
\newcommand{\R}{\mathbb{R}}
\newcommand{\Z}{\mathbb{Z}}
\newcommand{\N}{\mathbb{N}}
\newcommand{\p}{\partial }
\newcommand{\ee}{\mathrm{e}}
\newcommand{\uloc}{{\mathrm{uloc}}}
\newcommand{\iin}{{\mathrm{in}}}
\newcommand{\out}{{\mathrm{out}}}
\newcommand{\na}{\nabla}
\newcommand{\lec }{\lesssim }
\newcommand{\eqnb}{\begin{equation}}
\newcommand{\eqne}{\end{equation}}
\renewcommand{\d}{\mathrm{d}}
\newcommand{\LulocR}{L_{\uloc,R}}

\newcommand\Lu{L_{\uloc}}
\newcommand\Lut{L_{t,x-\uloc}}
\newcommand\Luloc{L_{3-\uloc}}

\newtheorem{thm}{Theorem}[section]
\newtheorem{lem}[thm]{Lemma}
\newtheorem{cor}[thm]{Corollary}
\newtheorem{prop}[thm]{Proposition}
\theoremstyle{remark}

\begin{document}
\title[Quantitative estimates for axisymmetric NSE]{Quantitative control of solutions to the axisymmetric Navier-Stokes equations in terms of the weak $L^3$ norm}
\author{ W. S. O\.za\'nski, S. Palasek}
\maketitle

\date{}

\medskip

\begin{abstract}
We are concerned with strong axisymmetric solutions to the $3$D incompressible Navier-Stokes equations. We show that if the weak $L^3$ norm of a strong solution $u$ on the time interval $[0,T]$ is bounded by $A \gg 1$ then for each $k\geq 0 $ there exists $C_k>1$ such that $\| D^k u (t) \|_{L^\infty (\mathbb{R}^3) } \leq t^{-(1+k)/2}\exp \exp A^{C_k}$ for all $t\in (0,T]$. 
\end{abstract}

\section{Introduction}
We are concerned with the 3D incompressible Navier-Stokes equations,
\eqnb\label{NSE}
\begin{cases}
&\p_t u - \Delta u + (u\cdot \nabla ) u + \nabla p =0,\\
&\div\,u=0 \qquad \text{ in } \R^3
\end{cases}
\eqne
for $t\in [0,T)$. While the question of global well-posedness of the equations remains open, it is well-known that the unique strong solution on a time interval $[0,T)$ can be continued past $T$ provided a regularity criterion holds, such as $\int_0^T \| \curl\, u \|_\infty dt <\infty $ (the Beale-Kato-Majda \cite{bkm} criterion), Lipschitz continuity up to $t=T$ of the direction of vorticity (the Constantin-Fefferman \cite{cf} criterion), or if $\int_0^T \| u \|_p^q dt< \infty$ for any $p\in [3,\infty ]$, $q\in [2,\infty ]$ such that $2/q+3/p\leq 1$ (the Ladyzhenskaya-Prodi-Serrin condition), among many others. The non-endpoint case $q<\infty $ of the latter condition was settled in the 1960s \cite{ladyzhenskaya_1967,serrin_1963,prodi_1959}, while the endpoint case $L^\infty_t L^3_x$ was only settled many years later by Escauriaza, Seregin, and \v{S}ver\'ak \cite{ess_2003}. The main difficulty of the endpoint case is related to the fact that $L^3$ is a critical space for $3$D Navier--Stokes, and \cite{ess_2003} settled it with an argument by contradiction using a blow-up procedure and new unique continuation results. This result implies that if $T_0>0$ is a putative blow-up time of \eqref{NSE}, then $\| u(t) \|_3$ must blow-up at least along a sequence of times  $t_k \to T_0^-$. While Seregin \cite{seregin_2011} showed that the $L^3$ norm must blow-up along any sequence of times converging to $T_0^-$, the question of quantitative control of the strong solution $u$ in terms of the $L^3$ norm remained open until the recent breakthrough work by Tao \cite{tao}, who showed that
\eqnb\label{tao_result}
| \na^j u(x,t) |\leq \exp \exp \exp (A^{O(1)}) t^{-\frac{j+1}2}
\eqne
for all $t\in [0,T]$, $j=0,1$, $x\in \R^3$, whenever
\[
\| u \|_{L^\infty ([0,T];L^3 (\R^3))} \leq A
\]
for some $A\gg 1$. This result implies in particular a lower bound
\[
\limsup_{t\to T_0^-} \frac{\| u(t) \|_3 }{\left( \log \log \log (T_0-t)^{-1}) \right)^c} = \infty,
\]
where $c>0$ and $T_0>0$ is the putative blow-up time,  and has subsequently been improved in some settings. For example, Barker and Prange \cite{barker_prange} and Barker \cite{barker} provided remarkable local quantitative estimates, and the second author \cite{palasek} proved that, in the case of axisymmetric solutions, 
\[
| \na^j u(x,t) |\leq \exp \exp (A^{O(1)}) t^{-\frac{j+1}2}
\]
for all $t\in [0,T]$, $j=0,1$, $x\in \R^3$, whenever
\[
\left\| r^{1-\frac3p }u \right\|_{L^\infty ([0,T];L^p (\R^3))} \leq A
\]
for some $A\gg 1$, $p\in (2,3]$.  In another work \cite{p2} he generalized \eqref{tao_result} to higher dimensions ($d\geq 4$), where, due to an issue related to the lack of Leray's intervals of regularity, one obtains an analogue of \eqref{tao_result} with four exponential functions. Recently Feng, He, and Wang \cite{feng} extended \eqref{tao_result} to the non-endpoint Lorentz spaces $L^{3,q}$ for $q<\infty$. We emphasize that all these generalizations rely on the same stacking argument by Tao \cite{tao}. In particular, the argument breaks down for the endpoint case $q=\infty$. 

\subsection{Tao's stacking argument and Type I blow-up}

In order to illustrate the issue at the endpoint space $L^{3,\infty}$, let us recall that the main strategy of Tao \cite{tao} is to show that if $u$ concentrates at a particular time, then there exists a widely separated sequence of length scales $(R_k)_{k=1}^K$ and $\alpha=\alpha(A)>0$ such that $\|u\|_{L^3(\{|x|\sim R_k\})}\geq\alpha$ for all $k$, which implies that
\eqnb\label{stacking}
\| u \|_3^3 = \int_{\R^3} |u|^3 \geq \sum_{k} \int_{|x|\sim R_k} |u|^3 \geq \alpha^3 K.
\eqne
The more singularly $u$ concentrates at the origin, the larger one can take $K$; thus the $L^3$ norm controls the regularity of $u$. More precisely, if $\|u\|_3\leq A$ and $u$ concentrates at a large frequency $N$ at time $T$, then one can take $\alpha=\exp(-\exp(A^{O(1)}))$ and $K\sim \log(NT^\frac12)$, which, by \eqref{stacking}, implies that $N\leq T^{-\frac12}\exp\exp\exp(A^{O(1)})$. This controls the solution in the sense that higher frequencies do not admit concentrations, and so a simple argument \cite[Section~6]{tao} implies the conclusion \eqref{tao_result}.\\

Let us contrast this $L^3$ situation with that of general Lorentz spaces $L^{3,q}$ with interpolation exponent $q\geq3$. In that case, $\|u\|_{L^{3,q}(\{|x|\sim R_k\})}\geq\alpha$ implies
\eqn{
\|u\|_{L^{3,q}(\mathbb R^3)}\gtrsim \big\|\|u\|_{L^{3,q}(\{|x|\sim R_k\})}\big\|_{\ell_k^q}\geq\alpha K^{\frac1q},
}
and so one should expect the bounds from the stacking argument \eqref{stacking} used in the Lorentz space $L^{3,q}$ extension \cite{feng} to degenerate as $q\to\infty$.
Indeed, if $|u(x)|=|x|^{-1}$  then, for some constant $\alpha >0$, we have $\| u \|_{L^{3,\infty } (\{|x|\sim R\})} \geq \alpha  $ for all $R>0$, yet $\| u \|_{L^{3,\infty } (\R^3)} \sim 1$ which shows that the first inequality in \eqref{stacking} fails for the $L^{3,\infty}$ norm. For this reason, the approach of Tao \cite{tao} (and, for related reasons, of Escauriaza-Seregin-\v Sver\'ak) to the $L^3$ problem cannot be extended to $L^{3,\infty}$.\\

This issue is in fact closely related to the study of Type 1 blow-ups and approximately self-similar solutions to \eqref{NSE}. Leray famously conjectured the existence of backwards self-similar solutions that blow up in finite time, a possibility later ruled out by Ne\v{c}as, R\r{u}\v{z}i\v{c}ka, and \v{S}ver\'{a}k \cite{necas_sverak} for finite-energy solutions and by Tsai \cite{tsai} for locally-finite energy solutions. The latter reference identifies the following as a very natural ansatz for blow-up:
\eq{\label{selfsimilar}
u(t,x)=\frac1{(T_0-t)^{\frac12}}U\left(\frac x{(T_0-t)^\frac12}\right),\quad U(y)=a\left(\frac y{|y|}\right)\frac1{|y|}+o\left(\frac1{|y|}\right)\text{ as }|y|\to\infty,
}
where $a:S^2\to\mathbb R^3$ is smooth. While Tsai \cite{tsai} shows that there are no solutions \emph{exactly} of this form, solutions that approximate this profile or attain it in a discretely self-similar way are promising candidates for singularity formation, as demonstrated by, for example, the Scheffer constructions \cite{o_scheffer1,o_scheffer2,scheffer_1985,scheffer_1987}, and the recent numerical evidence of an approximately self-similar singularity for the axisymmetric system due to Hou \cite{hou}. Unfortunately, criteria pertaining to $L^3$ such as those in \cite{ess_2003,tao,palasek} are less effective at controlling such solutions because $|x|^{-1}\notin L^3(\Rd)$, which shows the relevance of the weak norm $L^{3,\infty}$.\\

Specializing to the case of axial symmetry, it is known, for instance, that certain critical pointwise estimates of $u$ with respect to the distance from the axis imply regularity \cite{CSYT1,CSYT2,pan}. Moreover, Koch, Nadirashvili, Seregin, and \v{S}ver\'{a}k \cite{KNSS} proved a Liouville-type theorem for ancient axisymmetric solutions. Furthermore, Seregin \cite{seregin_2020} proved that finite-time blow-up cannot be of Type I. Thus, roughly speaking,  no axisymmetric solution can approximate the profile \eqref{selfsimilar} all the way up to a putative blow-up time $T_0$. However, this regularity is only qualitative (indeed, the proof uses an argument by contradiction based on a ``zooming in'' procedure), and so explicit bounds on the solution have not been available. \\

The main purpose of this work is to make this regularity quantitative, in a similar sense in which Tao \cite{tao} quantified the Escauriaza-Seregin-\v Sver\'ak theorem \cite{ess_2003}. This allows us to not only to rule out Type I singularies, but also to control how singular they can possibly become. For example, it lets us estimate the length scale up to which a solution can be approximated by a self-similar profile, see Corollary~\ref{cor_ss} for details. 
\subsection{The main regularity theorem}

We suppose that a strong solution to \eqref{NSE} on the time interval $[0,T]$ is axisymmetric, namely that
\begin{equation}\label{NSE_axisym}
\dd_\theta u_r=\dd_\theta u_3=\dd_\theta u_\theta=0,
\end{equation}
where $u_r, u_\theta, u_3$ denote (respectively) the radial, angular, and vertical components of $u$, so that
\[
u=u_r e_r + u_{\theta }e_\theta + u_3 e_3
\]
in cylindrical coordinates, where $e_r$, $e_\theta$, $e_3$ denote the cylindrical basis vectors. We assume further that $u$ remains bounded in $L^{3,\infty }$,
\eq{
\|u\|_{L^\infty ([0,T];L^{3,\infty}(\Rd))}\leq A\label{weakL3}
}
for some $A\gg 1$. We prove the following. 

\begin{thm}[Main result]\label{regularity}
Suppose $u$ is a classical axisymmetric solution of \eqref{NSE} on $[0,T]\times\Rd$ obeying \eqref{weakL3}. Then 
\eqn{
\|\grad^ju(t)\|_{L_x^\infty(\Rd)}\leq t^{-\frac{1+j}2}\exp\exp(A^{O_j(1)})
}
for all $j\geq 0$, $t\in [0,T]$.
\end{thm}

We note that, although our proof of the above theorem does use some of the basic a priori estimates (see Section~\ref{sec_basic_ests}) pointed out by Tao \cite{tao}, it follows a completely different scheme. Our main ingredients are parabolic methods applied to the swirl $\Theta \coloneqq r u_\theta $ near the axis, as well as localized energy estimates on 
\eqnb\label{phi_gamma_def}
\Phi \coloneqq \frac{\omega_r}{r} \quad \text{ and }\quad \Gamma \coloneqq \frac{\omega_\theta}r.
\eqne
In a sense, we use those estimates to  replace the Carleman inequalities appearing in Tao's \cite{tao} approach.

To be more precise, our proof builds on the work of Chen, Fang, and Zhang \cite{cfz_2017}, who showed that the energy norm of $\Phi$, $\Gamma$, 
\eqnb\label{phi_gamma_energy}
\| \Phi \|_{L^\infty_t L^2_x} + \| \Gamma \|_{L^\infty_t L^2_x} + \| \nabla \Phi \|_{L^2_t L^2_x} + \| \nabla \Gamma \|_{L^2_t L^2_x}  ,
\eqne
controls $u$ via an estimate on $\| u_\theta^2/r \|_{L^2}$ (see \cite[Lemma~3.1]{cfz_2017}). They also observed that one can indeed estimate this energy norm as long as the angular velocity $u_\theta $ remains small in any neighbourhood of the axis, namely if 
\eqnb\label{swirl_smallness_cfz}
\| r^d u_\theta \|_{L^\infty_t ([0,T]; L^{3/(1-d)}(\{ r \leq \alpha \} ))} \text{ is sufficiently small for some }\alpha >0\text{ and }d\in (0,1).
\eqne
In fact, this can be observed from the PDEs satisfied by $\Phi$, $\Gamma$, 
\eqnb\label{phi_gamma_pdes}
\begin{split}
&\left(  \p_t + u \cdot \na - \Delta - \frac2r \p_r \right) \Gamma + \frac{2}{r^2} u_\theta \omega_r =0,\\
&\left(  \p_t + u \cdot \na - \Delta - \frac2r \p_r \right) \Phi  - (\omega_r\p_r + \omega_3 \p_3 ) \frac{u_r}r=0,
\end{split}
\eqne
which show that, in order to control the energy of $\Gamma$, $\Phi$ one needs to control $u_r/r$, $\omega_r$, $\omega_3$ and $u_\theta$. However, $u_r/r$ can be controlled by $\Gamma$ in the sense that 
\eqnb\label{ur_trick}
\frac{u_r}r = \Delta^{-1} \p_3 \Gamma - 2 \frac{\p_r }{r} \Delta^{-2} \p_3 \Gamma
\eqne
(see \cite[p. 1929]{cfz_2017} for details), which is one of the main properties of function $\Gamma$. In particular, \eqref{ur_trick} lets us use the  Calder\'on-Zygmund inequality to obtain that
\eqnb\label{CZ_for_ur}
\left\| D^2 \frac{u_r}r \right\|_{L^q} \leq \| \p_3 \Gamma \|_{L^q}
\eqne
for $q\in (1,\infty)$ (see \cite[Lemma~2.3]{cfz_2017} for details). Moreover $\omega_r = r \Phi$, and $\omega_3 = \p_r (r u_\theta )/r$, which shows that the $L^2$ estimate of $\Phi$, $\Gamma$ relies only on control of $u_\theta$. In fact, away from from the axis, one can easily control $u_\theta$, while near the axis the smallness condition \eqref{swirl_smallness_cfz} is required in an absorption argument by the dissipative part of the energy, see \cite[(3.11)--(3.14)]{cfz_2017} for details.

In this work we obtain such control of $u_\theta$ thanks to the weak-$L^3$ bound \eqref{weakL3}, by utilizing parabolic theory developped by Nazarov and Ural'tseva \cite{nu_2012} in the spirit of the Harnack inequality. Namely, noting that the swirl $\Theta \coloneqq ru_\theta $ satisfies the autonomous PDE
\eq{\label{swirl_pde}
\Big(\dd_t+\Big(u+\frac2re_r\Big)\cdot\grad-\Delta\Big)\Theta=0
}
everywhere except for the axis, one can deduce (as observed in \cite[Section~4]{nu_2012}) H\"older continuity of $\Theta $ near the axis. A similar observation, but in a case of limited regularity of $u$ was used by Seregin \cite{seregin_2020} in his proof of no Type I blow-ups for axisymmetric solutions. We quantify this approach (see Proposition~\ref{prop_holder} below) to obtain an estimate on the H\"older exponent in terms of the weak-$L^3$ norm, and hence we obtain sufficient control of the swirl $\Theta$ in a very small neighbourhood of the axis. As for the outside of the neighbourhood, we obtain pointwise estimates on $u$ and all its derivatives, which are quantified with respect to $A$, and which improve the second author's estimates \cite[Proposition~8]{palasek}. This would enable one to close the energy estimates for the quantities in \eqref{phi_gamma_energy} if there exist sufficiently many starting times where the energy norms are finite. Indeed, given a weak $L^3$ bound \eqref{weakL3} and short time control of the dynamics of the energy \eqref{phi_gamma_energy}, control of $\| \Phi (T) \|_{L^2} + \| \Gamma (T) \|_{L^2}$ can be propagated from an initial time very close to $t=T$. Unfortunately, there are no times when we can explicitly control these energies in terms of $A$ due to lack of quantitative decay in the $x_3$ direction. The standard approach of propagating $L^2$ control of $\Phi,\Gamma$ from the initial data at $t=0$ (for instance, as in \cite{cfz_2017}) would lead to additional exponentials in Theorem~\ref{regularity}. \\

To avoid this issue and prove efficient bounds, we replace \eqref{phi_gamma_energy} with $L^2$ norms that measure $\Phi$ and $\Gamma$ uniformly-locally in $x_3$: namely, we consider
\eqnb\label{phi_gamma_energy_3uloc}
\| \Phi \|_{L^\infty_t L^2_{3-\uloc}} + \| \Gamma \|_{L^\infty_t L^2_{3-\uloc}} + \| \nabla \Phi \|_{L^2_t L^2_{3-\uloc}} + \| \nabla \Gamma \|_{L^2_t L^2_{3-\uloc}}, 
\eqne
where $\|\cdot \|_{L^2_{3-\uloc} }\coloneqq \sup_{z\in\mathbb R}\|\cdot \|_{L^2(\mathbb R^2\times[z-1,z+1])}$. See Proposition~\ref{prop_energy} below for an estimate of such energy norm. This approach gives rise to two further challenges. 

One of them is the $x_3$-$\uloc$ control of the solution $u$ itself in terms of \eqref{phi_gamma_energy_3uloc}. We address this difficulty by  an $x_3$-$\uloc$ generalization of the $L^4$ estimate on $u_\theta/ r^{1/2}$ introduced by \cite[Lemma~3.1]{cfz_2017}, together with a $x_3$-$\uloc$ bootstrapping via $\| u \|_{L^\infty_t L^6_{3-\uloc}}$, as well as an inductive argument for the norms $\| u \|_{L^\infty_t W_{\uloc}^{k-1,6}}$ with respect to $k\geq 1$, where ``$\uloc$'' refers to the uniformly locally integrable spaces (in all variables, not only $x_3$). We refer the reader to Steps 2--4 in Section~\ref{sec_pf_thm1} for details.

Another challenge is an $x_3$-$\uloc$ estimate on $u_r$ in terms of $\Gamma$. To be more precise, instead of the global estimate \eqref{CZ_for_ur}, we require $L^2_{3-\uloc}$   control of $u_r/r$, which is much more challenging, particularly considering the bilaplacian term in \eqref{ur_trick} above. To this end we develop a bilaplacian Poisson-type estimate in $L^2_{3-\uloc}$ (see Lemma~\ref{lem_doubleL}), which enables us to show   that
\eqnb\label{u_r_est_intro}
\left\|\nabla \p_r \frac{u_r}r\right\|_{L^2_{3-\uloc}}+ \Big\|\nabla \p_3 \frac{u_r}r\Big\|_{L^2_{3-\uloc}}\lesssim\|\Gamma\|_{L^2_{3-\uloc}}+\|\grad\Gamma\|_{L^2_{3-\uloc}},
\eqne
see Lemma~\ref{lem_localgamma}. Note that this is a $x_3$-$\uloc$ generalization of \eqref{CZ_for_ur}, and also requires the whole gradient on the right-hand side, rather than $\p_3\Gamma$ only. Such an estimate  lets us close the estimate of \eqref{phi_gamma_energy_3uloc}, and thus control all subcritical norms of $u$ in terms of $\| u \|_{L^{3,\infty }}$ (see Section~\ref{sec_pf_thm1} for details). \\

Having overcome the two difficulties of controlling the energy \eqref{phi_gamma_energy_3uloc}, we deduce (in \eqref{B}) that $\| \Gamma (t) \|_{ L^2_{3-\uloc }} \leq \exp \exp A^{O(1)}$ for all $t\in [1/2,1]$, whenever a solution $u$ satisfies $\|u\|_{L^\infty ([0,1];L^{3,\infty})}\leq A$; see Figure~\ref{fig1} (supposing that $T=1$). This suffices for iteratively improving the quantitative control of $u$ until $t=1$. Indeed, we first deduce a subcritical bound on the swirl-free part of the velocity on the same time interval, namely that $\| u_r e_r +u_ze_z\|_{L^p_{3-\uloc}} \lec_p \exp \exp A^{O(1)}$ for $p\geq 3$ and $t\in[1/2,1]$. We can then control (in \eqref{utheta/r12}) the time evolution of $\| u_\theta r^{-1/2} \|_{L^4_{3-\uloc}}$ over short time intervals, and so, choosing $t_0\in [0,1]$ sufficiently close to $1$ (by picking a time of regularity, see Lemma~\ref{lem_choice}) we then obtain 
(in \eqref{done}) that $\| u_\theta r^{-1/2} \|_{L^4_{3-\uloc}}$ and $\| u \|_{L^6_{3-\uloc }}$ are bounded by $\exp \exp A^{O(1)}$ for all $t\in [t_0,1]$, see Figure~\ref{fig1}. This subcritical bound allows one to also estimate $\| u \|_{W^{k,6}_{\uloc}}\leq \exp \exp A^{C_k}$ for every $k$, on a time interval of the same size (see Step~4 in Section~\ref{sec_pf_thm1}), which yields the claim of Theorem~\ref{regularity}.

\begin{center}
 \includegraphics[width=13cm]{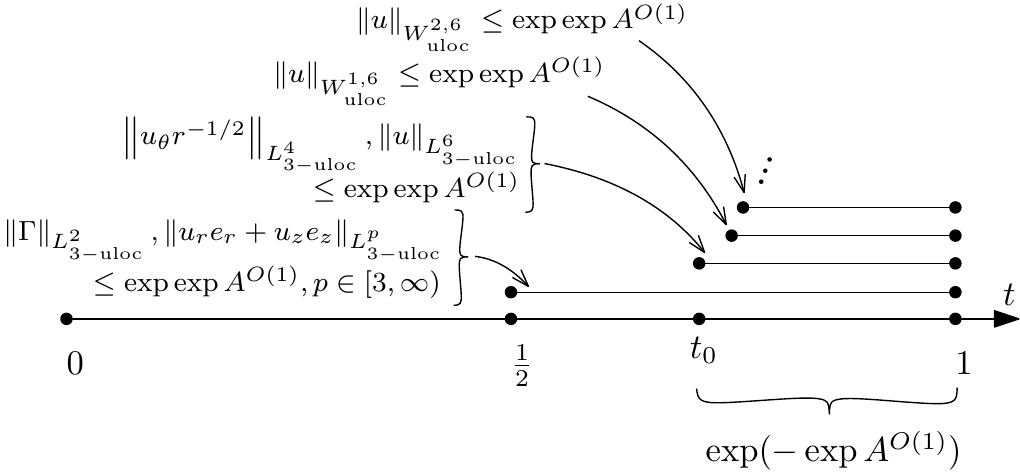}
 \end{center}
 \nopagebreak
 \captionsetup{width=.8\linewidth}
  \captionof{figure}{A sketch of the proof of Theorem~\ref{regularity}.}\label{fig1}

\subsection{A comparison of the blow-up rate}

We note that Theorem~\ref{regularity}, together with the well-known blow-up criterion $\| u(t) \|_\infty \geq c/(T_0-t)^{1/2}$ (see \cite[Corollary~6.25]{op}, for example), where $T_0>0$ is a putative blow-up time, immediately implies the following lower bound on the blow-up rate of $\| u(t) \|_{L^{3,\infty }}$.

\begin{cor}[Blow-up rate of the weak-$L^3$ norm]\label{corollary}
If $u$ is a classical axisymmetric solution of \eqref{NSE} that blows up at $T_0$, then
\eqnb\label{blowup_rate}
\limsup_{t\to T_0^-}\frac{\|u(t)\|_{L^{3,\infty}(\Rd)}}{(\log\log(T_0-t)^{-1})^c}=+\infty.
\eqne
\end{cor}

This corollary is also a consequence of a recent theorem of Chen, Tsai, and Zhang \cite{ctz}, who prove\footnote{Let us note the existence of a substantial misprint in the published version of \cite{ctz}: in their Theorem 1.4, as in our Corollary~\ref{corollary}, the blow-up rate is \emph{double}-logarithmic.}
\[
\limsup_{t\to T_0^-}\frac{\|b(t)\|_{\dot B^{-1}_{\infty, \infty }(\R^3)}}{\left(\log\log \frac{100}{T_0-t}\right)^{\frac1{48}-}}=+\infty,
\]
where $b\coloneqq u_r e_r + u_3 e_3$ denotes the swirl-less part of the velocity field $u$ (see \cite[Section 3.3]{lemarie_recent} for the relevant definition of $\dot B_{\infty,\infty}^{-1}$). Thus, since $\dot B^{-1}_{\infty, \infty }(\R^3)\supset L^{3,\infty}$, the above blow-up rate implies \eqref{blowup_rate}. We conjecture that a variant of Theorem~\ref{regularity} holds with the weak-$L^3$ norm replaced by such a critical Besov norm and can be proved using the ideas presented here.

In order to describe the relation of Corollary~\ref{corollary} to  \cite{ctz}, we note that the argument in \cite{ctz} proceeds by proving a pointwise estimate of the form 
\eqnb\label{ctz_claim}
| r u_\theta | \leq C \exp (-c |\log r |^\tau),
\eqne
where $c,C>0$, $\tau \in (0,1)$, for axisymmetric solutions obeying the slightly supercritical bound
\[
\frac{1}{R^\frac12}\| u \|_{L^\infty ((-R^2 ,0); L^2 (B_R) )}\leq K \left( \log \log \frac{100}R \right)^\beta\quad\text{ for all }R\in(0,1/4]
\]
for some $\beta\in(0,\frac18)$ and $K>0$. This is yet another application of Harnack inequality methods to axisymmetric Navier-Stokes equations. Rather than proving H\"older continuity of $\Theta $ under a global control of a critical norm as we do in Proposition~\ref{prop_holder}, \cite{ctz} obtains \eqref{ctz_claim} by an ``almost H\"older continuity,''
\eqnb\label{ctz_claim2}
\osc_{Q_\rho } \Theta \leq \exp \left(-c \left( \left( \log \frac{100}\rho \right)^\tau - \left( \log \frac{100}R \right)^\tau \right) \right) \osc_{Q_R} \Theta
\eqne
for $0<\rho < R\leq 1/4$, $\tau \in (0,1)$; see \cite[Proposition~1.2]{ctz}. A similar result in the case of $\tau=1/4$ has been obtained independently by Seregin \cite[Proposition~1.3]{seregin_2022}. Note that the case of $\tau=1$ corresponds to H\"older continuity.\\

We emphasize that the main point of our work is not to improve the blow-up rate but to give an explicit bound on $u$ and its derivatives in terms of only the critical norm---this is a strictly stronger result in the sense that it pertains to \emph{all} axisymmetric classical solutions, even those not blowing up. A na\"ive attempt to prove a similar quantitative theorem (e.g., using ideas of estimating axisymmetric vector fields from \cite{lei2017criticality}) would lead to a bound which, compared to Theorem~\ref{regularity}, would contain more iterated exponentials as well as severe dependence on the time $t$ and subcritical norms of the initial data. Instead, Theorem~\ref{regularity} parallels the results in \cite{tao} and improves on those in \cite{palasek} in the sense that the final bound depends only on $\|u\|_{L_t^\infty L_x^{3,\infty}}$ and a dimensional factor in $t$. This also leads to additional interesting corollaries: for instance, an explicit rate of convergence for $u(t)\to0$ as $t\to+\infty$, and the non-existence of nontrivial ancient axisymmetric solutions in $L_t^\infty L_x^{3,\infty}$. \\

A comparison of these results with the work of Chen, Tsai, and Zhang \cite{ctz} raises the following question: Is it possible to efficiently control (in the sense of Theorem~\ref{regularity}) $u$ and its derivatives in terms of only $b$ measured in some critical norm? In fact, in our proof of H\"older continuity of $\Theta $ near the axis (Proposition~\ref{prop_holder}) one can easily replace \eqref{weakL3} with boundedness of $\|b(t) \|_{L^{3,\infty}}$ in time, since ``$u$'' in \eqref{swirl_pde} can be replaced by ``$b$'', due to axisymmetry. However, we do require $L^{3,\infty}$ control of all components of $u$ for  other quantitative estimates leading to Theorem~\ref{regularity}. These include the basic estimates (Lemmas~\ref{lem_choice}--\ref{apriori}), quantitative decay away from the axis (Proposition~\ref{pointwisebounds}), as well as energy estimates on $\Gamma$ and $\Phi$ (Proposition~\ref{prop_energy}) and their implementation in the main argument (Section~\ref{sec_pf_thm1}).\\

A related open problem is to explicitly control $u$ in terms of $u_\theta$ only. In fact, despite a number of works \cite{cfz_2017,ladyzhenskaya_1968,lei_zhang_2011, np_2001,seregin_2022,ukhovskii_yudovich_1968} on the properties of the swirl $ru_\theta $, its role  in the regularity problem of axisymmetric solutions remains unclear.

\subsection{An estimate on the self-similar length scale}

One of the remarkable consequences of the quantitative estimate provided by Theorem~\ref{regularity} above is that it provides an estimate on the length scale up to which an axisymmetric solution to the NSE \eqref{NSE} can be approximated by a self-similar profile as in \eqref{selfsimilar}.

In order to make this precise, we will say that a vector field $b\in L^\infty(\Rd; \Rd)$ is \emph{nearly-spherical} if there exists $\delta \in (0,1/2)$  such that for every $R>0$, there exists $x_0\in\Rd$ with $|x_0|=R$ such that
\begin{align}\label{profilecondition}
|b(x_0)|\geq \frac{\|b\|_\infty}2\quad\text{and}\quad|b(x)-b(x_0)|\leq\frac{\|b\|_\infty}{4}\quad\text{for all }x\in B(x_0,\delta |x_0|).
\end{align}
Clearly any spherical profile $b(x)=a(x/|x|)$ is nearly-spherical for every $a\in C(\partial B(0,1))$ (in which case the choice of $\delta $ for  \eqref{profilecondition} to hold can be made  by a simple continuity argument). 
Let $\psi \in C_c^\infty (\Rd ; [0,1])$ be such that $\int \psi=1$, and let $\psi_l (x)\coloneqq l^{-3} \psi (x/l)$ denote a mollifier at a given length scale $l>0$. We also set $\widetilde{\psi_l} \coloneqq \psi_l \ast \psi_l$.

We note that, letting $R\coloneqq 2l/\delta $, we can find   $x_0\in \Rd$ with  $|x_0|= 2l/\delta $ and satisfying \eqref{profilecondition}. In particular
\[
\left| \left( \widetilde{\psi_l} \ast \frac{b(\cdot )}{|\cdot |}\right) (x_0) \right| =\left|\int_{B(x_0,2l )}\widetilde{\psi_l}(x_0-y)\frac{b(y)}{|y|}\d y\right|\gtrsim \frac{|b(x_0)| - \| b \|_\infty /4}{(1+\delta )|x_0|} \geq \frac{\delta \| b \|_\infty }{16 l},
\]
which shows that
\eqnb\label{linfty_lowerbd}
\left\| \widetilde{\psi_l } \ast \frac{b(\cdot )}{|\cdot |} \right\|_\infty \geq \frac{\delta \| b \|_\infty }{16 l} 
\eqne
for every length scale $l>0$.
This simple fact lets us deduce from Theorem~\ref{regularity} that, if an axisymmetric solution approximates a self-similar profile $b(t,x)/|x|$ up to length scale $l(t)$, where $b$ is nearly-spherical uniformly on $[0,t]$, then $l(t)$ cannot be smaller than a particular quantitative threshold.  
\begin{cor}\label{cor_ss}
If $u$ is a strong axisymmetric solution $u$ of \eqref{NSE} on $[0,T]$,  
\eq{\label{approximateprofile}
\Big\|u(t )-\psi_{l(t )}*\frac{b (t,x)}{|x|}\Big\|_{L^{3,\infty}}\leq \sigma  \|b(t)\|_\infty
}
for $t\in[0,T]$, and $\sigma < c \delta $, where $c>0$ is a sufficiently small constant and $b(T)$ is nearly-spherical with constant $\delta$, then 
\eqn{
l (T)\gtrsim \delta T^{\frac12}\|b(T)\|_\infty\exp \left( -\exp \left( \|b\|_{L_{t,x}^\infty([0,T]\times\Rd)}^{O(1)}\right) \right).
}
\end{cor}

\begin{proof} 
We note that, at time $T$,
\eqn{
\|u\|_\infty&\gtrsim \|\psi_{l}*u\|_\infty\\
&\geq\left\|\widetilde{\psi}_{l}*\frac{b(\cdot )}{|\cdot |}\right\|_\infty-\left\|\psi_{l}*\left(u-\psi_l*\frac{b(\cdot )}{|\cdot |}\right)\right\|_\infty\\
&\geq \frac{\delta \| b \|_\infty }{16l } - C l^{-1}\left\| u-\psi_l*\frac{b(\cdot )}{|\cdot |} \right\|_{L^{3,\infty}}\\
&\geq \left( \frac{\delta }{16} - C\, \sigma \right) \frac{ \left\| b \right\|_\infty }{l}.
}
Thus $\| u (T)\|_\infty \geq \delta \| b(T) \|_\infty / 32l$ if $\sigma \in (0,\delta /32C ) $. Since also
\[
\| u (t)\|_{L^{3,\infty } } \leq \left\|\widetilde{\psi}_{l(t)}*\frac{b(t, \cdot )}{|\cdot |}\right\|_{L^{3,\infty } }+\left\|u(t)-\psi_{l(t)}*\frac{b(t,\cdot )}{|\cdot |}\right\|_{L^{3,\infty } } \leq C \| b(t,\cdot ) \|_{\infty }
\]
for \emph{all} $t\in[0,T]$, Theorem~\ref{regularity} implies that
\[
\frac{\delta \| b (T) \|_\infty }{32 \,l(T)} \leq \| u (T) \|_\infty \lec T^{-1/2} \exp \exp \left(  \| b \|_{L^\infty  ([0,T]\times \Rd)}^{O(1)}\right) , 
\]
from which the claim follows. 
\end{proof}

\subsection{Organization of the paper}

The structure of the paper is as follows. In the following Section~\ref{sec_prelims} we discuss preliminary concepts related to the  Lorentz spaces $L^{p,q}$, the Bogovski\u i operator, a simple Poisson-type tail estimate that we will later (in Section~\ref{sec_poisson}) expand to obtain our Poisson-type estimate \eqref{u_r_est_intro} above, as well as some properties of cylindrical coordinates. In Section~\ref{sec_axisym_fcns} we discuss some  properties of axisymmetric functions, including an axisymmetric Bernstein inequality (Section~\ref{sec_bernstein}) and a quantified version of Hardy's inequality (Section~\ref{sec_quantified_hardy}). In Section~\ref{sec_ests_nse} we present some quantitative estimates of the 3D Navier--Stokes equations, including the Picard iterates (Section~\ref{sec_picard}), times of regularity, bounded total speed, and second derivatives estimates (Section~\ref{sec_basic_ests}), all of which remain valid without the assumption of axisymmetry. The following section, Section~\ref{sec_ests_axisymm_nse}, is dedicated to quantitative estimates that are specific to the axisymmetric setting \eqref{NSE_axisym} of the equations \eqref{NSE}. These include the statement of the H\"older estimate of the swirl $\Theta$ mentioned above (Section~\ref{sec_holder}), pointwise estimates away from the axis (Section~\ref{sec_pointwise}), as well as the Poisson-type $x_3$-$\uloc$ estimate on $u_r/r$ \eqref{u_r_est_intro} (Section~\ref{sec_poisson}). In Section~\ref{sec_en_ests_cfz} we prove the energy estimate \eqref{phi_gamma_energy_3uloc} for $\Gamma$ and $\Phi$  mentioned above, and Section~\ref{sec_pf_thm1} combines the developed methods to prove the main theorem, Theorem~\ref{regularity}. Finally, Appendix~\ref{sec_NU} includes a detailed verification of the H\"older estimate of $\Theta$.

\section{Preliminaries}\label{sec_prelims}

Given $f\colon \Omega \to \R$ we let
\[
\osc_{\Omega }\, f \coloneqq \sup_{\Omega } f - \inf_{\Omega } f 
\]
denote the oscillation of $f$ over $\Omega$. We also denote by $\fint_{\Omega } \coloneqq \frac{1}{|\Omega |} \int_{\Omega }$ the average over $\Omega$.\\

We use standard definitions of Lebesgue spaces $L^p (\Omega ) $, Sobolev spaces $W^{k,p}(\Omega ) $, spaces of continuous functions $C(\Omega ) $, spaces $C_c(\Omega ) $ of continuous functions with compact support. For brevity of notation we often omit ``$\Omega$'' in the notation if $\Omega = \R^3$; for example $W^{1,\infty } \equiv W^{1,\infty } (\R^3)$. We use the convention $\| \cdot \|_p \coloneqq \| \cdot \|_{L^p (\R^3 )}$, and we reserve the notation $\| \cdot \| \coloneqq \| \cdot \|_{2}$ for the $L^2(\R^3)$ norm. We also write $\int \coloneqq \int_{\R^3}$. Given $p\in[1,\infty]$, we  define the uniformly local $L^p$ norms,
\eqnb\label{uloc_norms}
 \|u\|_{\Lu^p}\coloneqq \sup_{x\in\Rd}\|u\|_{L_x^p(B(x,1))}\quad \text{ and }\quad \|u\|_{\Lut^p}\coloneqq \big\|\|u\|_{\Lu^p}\big\|_{L_t^p},
 \eqne
as well as the norms that are uniformly local in $x_3$ only,
\eqnb\label{l2_3uloc_def}
\|f\|_{L^p_{3-\uloc} (\Rd)}\coloneqq \sup_{z\in\mathbb R}\|f\|_{L^p(\mathbb R^2\times[z-1,z+1])}.
\eqne

We let $\Psi (x,t) \coloneqq(4\pi t)^{-3/2} \ee^{-x^2/4t}$ denote the heat kernel, which satisfies
\eqnb\label{heat_bounds}
\| \na^k \Psi (t) \|_p = C_{k,p} t^{-\frac32 \left(1-\frac1p  \right) -\frac{k}{2} }.
\eqne
We often use the notation $\ee^{t\Delta }f \coloneqq \Psi (t) * f$.

Given $N\in \{ 2^k \colon k\in \mathbb{N}\}$ we let $P_N$ denote the $N$-th Littlewood-Paley projection. We recall a localized version of the Bernstein inequality
\eqnb\label{bernstein_localized}
\| P_N f \|_{L^q (\Omega ) } \lec_k N^{\frac{3}{p_1}-\frac3q } \| P_N f \|_{L^{p_1} (\Omega_{R})} + (RN)^{-k} N^{\frac{3}{p_2}-\frac3q } \| P_N f \|_{L^{p_2}},
\eqne
where $\Omega\subset \R^3$ is an open set, $k\geq 1$, $\Omega_{R } \coloneqq \{ x\in \R^3 \colon \mathrm{dist}(x,\Omega ) < R  \}$, $q\in [1,\infty ]$ and $p_1,p_2 \in [1,q]$; see \cite[Lemma~2.1]{tao} for a proof.

\subsection{Lorentz spaces}

We recall the Lorentz spaces, defined by  
\eqnb\label{lorentz_def}
\| f \|_{L^{p,q}} \coloneqq p^{1/q} \| \lambda | \{ |f |\geq \lambda  \} |^{1/p} \|_{L^q (\R_+ , \frac{\d \lambda }{\lambda } )}
\eqne
for $q<\infty$ and
\[
\| f \|_{L^{p,\infty}} \coloneqq  \| \lambda | \{ |f |\geq \lambda  \} |^{1/p} \|_{L^\infty (\R_+ , \frac{\d \lambda }{\lambda } )}.
\]

We recall the H\"older inequality for Lorentz spaces,
\eqnb\label{holder_lorentz}
\| fg \|_{L^{p,q}} \leq C_{p_1,p_2,q_1,q_2} \| f \|_{L^{p_1,q_1}} \| g \|_{L^{p_2,q_2}}, 
\eqne
whenever $1/p=1/p_1+1/p_2$, $1/q=1/q_1+1/q_2$, $p_1,p_2,p\in (0,\infty )$, $q_1,q_2,q\in (0,\infty ]$. We refer the reader to \cite[Theorem~6.9]{terry_notes} for a proof of \eqref{holder_lorentz}. The H\"older inequality can be very useful when estimating some localized integrals in terms of the $L^{p,\infty}$ norm. For example, if $\phi \in C_0^\infty (\Omega )$ is a smooth cutoff function then we have the simple estimate
\[
\| \phi \|_{L^{p,1}} = p\int_0^\infty | \{ |\phi |\geq \lambda \} |^{1/p} \d \lambda \leq p\int_0^{\| \phi \|_\infty } | \{ |\phi |\geq \lambda \} |^{1/p} \d \lambda \leq p |\Omega |^{1/p} \| \phi \|_\infty , 
\]
which shows that, for example
\[
\int_{\Omega } fg \leq \| f \|_{L^{3,\infty}} \| g \|_2 |\Omega |^{1/6}.
\]
 This simple method allows us to use the weak $L^3$ space to estimate some integrals over a region close to the axis of symmetry.
 
 We also note two Young's inequalities involving weak $L^p$ spaces
 \eqnb\label{young_weak}
 \| f\ast g \|_{L^{p,\infty }} \lec \| f \|_1 \| g \|_{L^{p,\infty }}\qquad \text{ for } p\in (1,\infty ),
 \eqne
  \eqnb\label{young_weak1}
 \| f\ast g \|_{p} \lec \| f \|_r \| g \|_{L^{q,\infty }} \quad \text{ for } p,q,r\in (1,\infty ) \text{ with } \frac1p +1 = \frac1q+\frac1r ,
 \eqne
 see \cite[Proposition~2.4(a)]{lemarie_recent} and \cite[Theorem~A.16]{nse_book} for details (respectively). 

\subsection{The Bogovski\u{\i} operator}

We recall that, given $p\in (1,\infty)$, an open ball $B\subset \R^3$,  $b\in W^{1,p}(B)$ such that $\div\, b =0$, and $\phi \in C_0^\infty (B; [0,1])$ such that $\phi=1 $ on $B/2$ there exists $\overline{b} \in W^{1,p} (\R^3)$ such that $\overline{b} =0$ outside $B$ and inside $B/2$, 
\eqnb\label{bogovskii}
\div \, \overline{b} = \div (\phi b) \quad \text{ and } \quad \| \overline{b } \|_{W^{1,p}} \lec_B \| b \|_{W^{1,p}(B)},
\eqne
due to the Bogovski\u{\i} lemma (see \cite{bogovskii_79,bogovskii_80} or \cite[Lemma~III.3.1]{galdi_book}, for example). Here we use the non-homogeneous $W^{1,p}$ norm and so the implicit constant in \eqref{bogovskii} may depend of the size of $B$. We note that the Bogovski\u{\i} lemma often assumes that the domain is star-shaped (which is not the case for $B\setminus B/2$), but it can be overcome in this particular setting by applying the partition of identity to $\phi$; see \cite[Section~2.3]{ozanski_2021} for example.

\subsection{A Poisson-type tail estimate}

Here we are concerned with a Poisson equation of the form $-\Delta f = D^2 g$, and we show that any $W^{k,\infty} (B(0,1))$ norm of $\na f$ can be bounded by the $L^1_{\uloc}$ norm of $g$, if $g=0$ on $B(0,2)$.

To be more precise, we let $\psi \in C_c^\infty (B(0,1);[0,1])$ be such that $\psi =1 $ on $B(0,1/2)$. Given $y\in \R^3$ we set
\eqnb\label{translate_notation}
\psi_y (x) \coloneqq \psi (x-y).
\eqne
and
\[
\widetilde{\psi } \coloneqq \sum_{\substack { j\in \Z^3 \\ |j|\leq 10 } } \psi_j.
\]
\begin{lem}\label{lem_poisson_tail}
Suppose that $f = D^2 (-\Delta )^{-1} (g (1-\widetilde{ \psi} ))$ for some $g\in L^2$. Then
\[
\| \psi \na f \|_{W^{k,\infty}} \lec_k \| g \|_{L^1_{\uloc}}  \qquad \text{ for } k\geq 0.
\]
\end{lem}
\begin{proof}
We note that
\[
\p_i f (x) = \int \frac{(x_i-y_i)g(y)(1-\tilde{\phi }(y))}{|x-y|^5} \d y
\]
for $x\in \supp\,\phi$, and so 
\[\begin{split}
|\na f (x) |&\leq \int_{\{ |x-y | \geq 5  \} }\frac{|g(y) |}{|x-y|^4}\d y \\
&\leq  \sum_{\substack{j\in \Z^3\\ |j|\geq 2 }} \int_{x_1+j_1}^{x_1+j_1+1} \int_{x_2+j_2}^{x_2+j_2+1} \int_{x_3+j_3}^{x_3+j_3+1} \frac{|g(y)|}{|x-y|^4} \d y_3\,\d y_2 \, \d y_1\\
&\lec \| g \|_{L^1_{\uloc} }\sum_{\substack{j\in \Z^3\\ |j|\geq 2 }} |j|^{-4} \lec \| g \|_{L^1_{\uloc} },
\end{split}
\]
as required. An analogous argument applies to higher derivatives of $f$.
\end{proof}
The above proof demonstrates a simple method of tail estimation which we will later use to obtain a $L^2_{3-\uloc}$ estimate of $u_r/r$ in terms of $\Gamma$, mentioned in the introduction (recall \eqref{u_r_est_intro}). In fact, to this end, a similar strategy can be applied in the $x_3$ direction only, and can be extended to the more challenging biLaplacian Poisson equation (see Lemma~\ref{lem_doubleL} below).  

\subsection{Cylindrical coordinates}\label{sec_cylindrical}

Given $x\in \R^3$ we denote by $x'\coloneqq (x_1,x_2)$ the horizontal variables, and $r\coloneqq (x_1^2 + x_2^2)^{1/2}$ denotes the radius in the cylindrical coordinates. We often use the notation
\[
\{ r < r_0 \} \coloneqq \{ x\in \R^3 \colon r<r_0 \}
\]
for a given $r_0>0$.

We recall a version of the Hardy inequality
\eqnb\label{hardy1}
\| r^{-1} f \|_{L^{q} (\Omega )} \lec C(\Omega ) \| f \|_{L^{q} (\Omega )} +\| \nabla f \|_{L^{q} (\Omega )},
\eqne
where $\Omega $ is a bounded domain and $q\in (1,2]$; see \cite[Lemma~2.4]{cfz_2017} for a proof.

We recall the divergence operator in cylindrical coordinates: if $v= v_r e_r + v_\theta e_\theta + v_3 e_3$ then
\eqnb\label{div_cylindr}
\div \, v = \frac{1}r \p_r (rv_r )+\frac 1r \p_\theta v_\theta + \p_3 v_3.
\eqne 

We say that a vector field $v$ is axisymmetric if \eqref{NSE_axisym} holds. In such case we have
\eqnb\label{LR_formula}
|\grad'v|^2=(\dd_r v_r)^2+(\dd_rv_\theta)^2+(\dd_rv_3)^2+\frac1{r^2}(v_r^2+v_\theta^2),
\eqne
which implies the pointwise bounds
\eqn{
\frac{|v_r|}{r},\,\frac{|v_\theta|}r\leq|\grad'v|.
}
Here $\nabla'$ refers to the gradient with respect to the horizontal variables $x'$ only.

Moreover,
\eqnb\label{drr_main}
|\p_{rr} f | \lec | D^2 f|.
\eqne
Indeed, since 
\[
\p_r = \cos \theta \, \p_1  +\sin \theta \, \p_2  = \frac{x_1}{|x'|} \p_1  + \frac{x_2}{|x'| } \p_2 ,
\]
where $x' \coloneqq  (x_1,x_2)$ refers to the horizontal variables, we can compute that  
\[
\begin{split}
\p_{rr} &= \frac{x_1^2}{|x'|^2} \p_{11} + 2 \frac{x_1 x_2}{|x'|^2 } \p_1 \p_2  + \frac{x_2^2}{|x'|^2 } \p_{22} ,
\end{split}
\]
from which \eqref{drr_main} follows.  More generally,
\[\begin{split}
\p_{rrr} &= \frac{x_1^3}{|x'|^3} \p_{111} + \frac{3x_1^2 x_2}{|x'|^3} \p_{11}\p_2 + \frac{3x_1x_2^2}{|x'|^3 } \p_1 \p_{22} + \frac{x_2^3}{|x'|^3} \p_{222} ,\\
\p_{rrrr} &= \frac{x_1^4}{|x'|^4} \p_{1111} + \frac{4x_1^3 x_2}{|x'|^4} \p_{111}\p_2 + \frac{6x_1^2x_2^2}{|x'|^4}\p_{11}\p_{22} + \frac{4x_1x_2^3}{|x'|^4 } \p_1 \p_{222} + \frac{x_2^4}{|x'|^4} \p_{2222}.
\end{split}
\]
This shows that
\eqnb\label{d4_in_r_x_3}
| D^3_{r,x_3} f | \lec | D^3 f|\qquad \text{ and } \qquad | D^4_{r,x_3} f | \lec | D^4 f|
\eqne
for any axisymmetric $f$ (here, for example, $D^4$ refers to all fourth order derivatives with respect to $x_1,x_2,x_3$).

\section{Properties of axisymmetric functions}\label{sec_axisym_fcns}

Here we discuss some properties of axisymmetric functions, including an axisymmetric Bernstein inequality and a quantified Hardy's inequality.

\subsection{Bernstein inequalities}\label{sec_bernstein}

Here we discuss a version of the axisymmetric Bernstein inequality provided by \cite[Proposition~1]{palasek} that involves the weak $L^3$ space.

\begin{lem}\label{bernstein}
Let $T_m$ be a Fourier multiplier whose symbol $m$ is supported on $B(0,N)$ with $|\grad^jm|\leq MN^{-j}$ and $1<q< p\leq\infty$. If either $-\frac2p<\alpha<\frac1q-\frac1p$ or $p=\infty$ and $\alpha=0$, we have
\eqn{
\|r^\alpha T_mu\|_{L^p}\lesssim MN^{\frac3q-\frac3p-\alpha}\|u\|_{L^{q,\infty}}
}
for all axisymmetric scalar- or vector-valued functions $u$.
\end{lem}

\begin{proof}
We normalize $M=N=1$. Under these assumptions on $p,\alpha$, Proposition 1 in \cite{palasek} implies
\eqn{
\|r^\alpha T_mu\|_{L^p}\lesssim\|P_{\leq10}u\|_{L^{q+\epsilon}}
}
since $T_mP_{\leq10}=T_m$, for an $\epsilon>0$ sufficiently small depending on $p,q,\alpha$. Let $\psi$ be the kernel such that $P_{\leq10}=\psi*$. Then by the weak Young inequality \eqref{young_weak1},
\eqn{
\|P_{\leq10}u\|_{L^{q+\epsilon}}\lesssim\|\psi\|_{L^{1+O(\epsilon)}}\|u\|_{L^{q,\infty}}\lesssim\|u\|_{L^{q,\infty}}.
}
\end{proof}

A useful consequence of the above lemma is the following heat kernel estimate
\eq{
\|r^\alpha \ee^{\Delta} \grad^jf\|_{L^p}&\leq\|r^\alpha \ee^{\Delta}\grad^jP_{\leq1}f\|_{L^p}+\sum_{N>1}\|r^\alpha \ee^{\Delta}\grad^jP_Nf\|_{L^p}\nonumber\\
&\lesssim_{\alpha, p ,q,j}\|f\|_{L^{q,\infty}}(1+\sum_{N>1}\ee^{-N^2/100}N^{j+\frac3q-\frac3p})\nonumber\\
&\lesssim_{ p ,q,j} \|f\|_{L^{q,\infty}}\label{heat}
}
under the same assumptions on the parameters as in Lemma \ref{bernstein}. 

\subsection{A quantified version of the Hardy inequality}\label{sec_quantified_hardy}

By the classical Hardy inequality 
\[
\| r^{-\frac3p+\frac12}f \|_p \lec_{p}  \left( \|f\|_{2}+\|\grad f\|_{2}\right)
\]
for any axisymmetric $f$, and $p\in (2,6)$ (see \cite[Lemma~2.6]{cfz_2017}, for example). Here we prove a version of this inequality, which is localized in the horizontal variables, ``uloc'' in $x_3$, and which has a quantified divergence of the constant near $p=2$. Namely we prove the following.

\begin{lem}[Quantified Hardy inequality]\label{lem_hardy}
For $p\in(2,6-\epsilon)$,
\eqn{
\|r^{-\frac3p+\frac12}f\|_{\Luloc^p(\{r\leq1\})}&\lesssim_\epsilon(p-2)^{-O(1)}\left( \|f\|_{\Luloc^2(\{r\leq 1\})}+\|\grad f\|_{\Luloc^2(\{r\leq 1\})}\right).
}
\end{lem}

\begin{proof}
From the Sobolev embedding
\eqn{
\|u\|_{L^{2p/(2-p)}(\mathbb R^2)}&\lesssim(2-p)^{-O(1)}\|\grad u\|_{L^p(\mathbb R^2)}
}
for $p<2$, (see, e.g., \cite{talenti} where the sharp constant is computed), one can prove the two-dimensional Gagliardo-Nirenberg inequality
\eq{\label{GN}
\|f\|_{L^q(B(1))}&\lesssim q\left( \|f\|_{L^6(B(1))}^{\frac6q}\|\grad f\|_{L^2(B(1))}^{1-\frac6q}+\|f\|_{L^p(B(1))}\right)
}
for $q>6$. Fix $\epsilon>0$ to be specified. Then
\eqn{
\Big\|\frac {f(\cdot , x_3)}{r^{\frac3q-\frac12}}\Big\|_{L_{x'}^q(r\geq\epsilon)}&\leq\|r^{-\frac3q+\frac12}\|_{L_{x'}^{6q/(6-q)}(\{r\geq\epsilon\})}\|f(\cdot , x_3)\|_{L_{x'}^6(\mathbb R^2)}\lesssim\epsilon^{-\frac1q+\frac16}\|f(\cdot , x_3)\|_{L_{x'}^6(\mathbb R^2)}.
}
Inside, for any $\frac1s\in(\frac3{2p}-\frac14,\frac1p)$, by \eqref{GN},
\eqn{
\Big\|\frac {f(\cdot , x_3)}{r^{\frac3p-\frac12}}\Big\|_{L_{x'}^p(r\leq\min(1,\epsilon))}&\leq\|r^{-\frac3p+\frac12}\|_{L_{x'}^s(r<\min(1,\epsilon))}\|f(\cdot , x_3)\|_{L_{x'}^{ps/(s-p)}(B(1))}\\
&\lesssim\Big(\frac1s-\frac3{2p}+\frac14\Big)^{-\frac1s}\Big(\frac1p-\frac1s\Big)^{-1}\\
&\quad\times\left( \epsilon^{-\frac3p+\frac12+\frac2s}\|f(\cdot , x_3)\|_{L_{x'}^6(B(1))}^{\frac6p-\frac6s}\|\grad f(\cdot , x_3)\|_{L_{x'}^2(B(1))}^{1-\frac6p+\frac6s}+\|f(\cdot , x_3)\|_{L_{x'}^p(B(1))}\right).
}
Upon taking $\epsilon=\|f\|_6^3/\|\grad f\|_2^3$ and $\frac1s=\frac4{3p}-\frac16$,
\eqn{
\Big\|\frac {f(\cdot , x_3)}{r^{\frac3p-\frac12}}\Big\|_{L_{x'}^p(B(1))}&\lesssim(p-2)^{-O(1)}\left( \|f(\cdot , x_3)\|_{L_{x'}^6(B(1))}^{\frac32-\frac3p}\|\grad f(\cdot , x_3)\|_{L_{x'}^2(B(1))}^{-\frac12+\frac3p}+\|f(\cdot , x_3)\|_{L_{x'}^p(B(1))}\right).
}
Finally by H\"older's inequality, Sobolev embedding, and Gagliardo-Nirenberg interpolation, we find
\eqn{
\Big\|\frac f{r^{\frac3p-\frac12}}\Big\|_{L_x^p(B_{\mathbb R^2}(1)\times B_{\mathbb R}(z,1))}&\lesssim(p-2)^{-O(1)}\|f\|_{H_x^1(B_{\mathbb R^2}(1)\times B_{\mathbb R}(z,1))},
}
as required.
\end{proof}

\section{Basic estimates for the Navier-Stokes solutions}\label{sec_ests_nse}

Here we discuss some estimates for the Navier-Stokes equations without the assumption of axisymmetry. 

\subsection{The Picard estimates}\label{sec_picard}
We define the flat and sharp Picard iterates
\eqnb\label{picard_def}
\ulin_n(t)\coloneqq \ee^{(t-t_n)\Delta}u(t_n)-\int_{t_n}^t \ee^{(t-t')\Delta}\mathbb P\div ( \ulin_{n-1}\otimes\ulin_{n-1}(t'))\d t',\quad\unlin_n \coloneqq u-\ulin_n
\eqne
for all $n=1,2,\ldots$ and $t\geq t_n$, where $t_n\in[0,\frac12)$ is an increasing sequence of times, and $\ulin_0\coloneqq 0$, $\unlin_0\coloneqq u$. We have the following.

\begin{lem}[Basic Picard estimates]\label{sharpflat}
Assume $u$ solves \eqref{NSE} on $[0,1]\times\Rd$ with the bound \eqref{weakL3}. If  $p\in (3, \infty ]$ and $-\frac2p<\alpha<\frac13-\frac1p$ or $p=\infty$ and $\alpha=0$, we have
\eq{\label{flatbounds}
\|r^\alpha\grad^j\ulin_n\|_{L_t^\infty L_x^p([\frac12,1]\times\Rd)}&\leq A^{O_{n,j,p}(1)},\\
\| \unlin_n\|_{L_{t}^\infty L_x^q([\frac12,1]\times\Rd)}&\leq A^{O_{n,q}(1)} \qquad \text{ for all }q\in (1,3), \label{sharpbound} \\
\|\grad^jP_N\ulin_n\|_{L_{t,x}^\infty([\frac12,1]\times\Rd)}&\leq \ee^{-N^2/O_{n,j}(1)} A^{O_{n,j}(1)},\label{frequencylocalized}
}
as well as the energy estimate
\eq{\label{energy}
\|\grad\unlin_n\|_{L_{t,x}^2([\frac12,1]\times\Rd)}\leq A^{O_n(1)}.
}
In particular,
\eq{\label{localenergy}
\|\grad u\|_{\Lut^2([\frac12,1]\times\Rd)}\leq A^{O(1)}.
}
\end{lem}
The proof of \eqref{flatbounds}--\eqref{frequencylocalized} above relies only on the definition \eqref{picard_def} as well as basic heat estimates \eqref{heat_bounds}, which, together with the weak Young's inequality \eqref{young_weak1}, can be used in the same way as \cite[(3.11)--(3.13)]{tao} and  \cite[Proposition~2.5]{p2} to obtain the estimates with  $\| u \|_{L^\infty ([0,1];L^{3,\infty })}\leq A$ on the right-hand side.

\subsection{Basic estimates}\label{sec_basic_ests}

Here we assume that $u$ satisfies \eqref{NSE} with the weak $L^{3,\infty}$ bound \eqref{weakL3} on the time interval $[0,T]$.

\begin{lem}[Choice of time of regularity]\label{lem_choice}
If $u$ solves \eqref{NSE} on a time interval $I$ and satisfies $\|u\|_{L_t^\infty L_x^{3,\infty}(I\times\Rd)}\leq A$, then there exists $t_*\in I$ such that
\eqn{
\|\grad^ju(t_*)\|_{L_{x}^\infty(\Rd)}&\leq |I|^{-\frac{1+j}2}A^{O(1)}
}
for all $j=0,1,2,\ldots,10$.
\end{lem}
\begin{lem}[Bounded total speed]
We have the bounded total speed estimate
\eqn{
\|u\|_{L_t^1L_x^\infty(I/2\times\Rd)}\leq |I|^{\frac12}A^{O(1)}.
}
\end{lem}
The 2 lemmas above follow by the same arguments in \cite[Lemma~3.1]{tao} and \cite[Propositions~3.1--2]{feng} using the estimates in Lemma \ref{sharpflat}. In particular, it is straightforward to check that the proofs of Propositions 3.1 and 3.2 in \cite{feng} are still valid in Lorentz spaces $L^{p,q}$ with $q =\infty$. Furthermore, we estimate $\nabla^2 u$ in terms of $A$.

\begin{lem}[2nd order derivatives estimates]\label{apriori}
If $u$ solves \eqref{NSE} on $[0,T]$ and obeys \eqref{weakL3}, then
\eqn{
\|\grad^2u\|_{\Lut^p([\frac T2,T]\times\Rd)}\lesssim_p A^{O(1)}T^{\frac5{2p}-\frac32}
}
for $p\in[1,\frac43)$, where the ``$\uloc$'' norm is considered as the supremum of the $L^p$ norms over $B(T^{1/2})\subset \R^3$ (instead of $B(1)$, recall \eqref{uloc_norms}).
\end{lem}

\begin{proof}
We use an approach due to Constantin \cite{constantin}. First rescale to make $T=1$. For every $\epsilon\in(0,\frac12)$, we define the approximation to the function $\langle x\rangle:=(1+|x|^2)^\frac12$,
\eqn{
q_\epsilon (x)\coloneqq \langle x\rangle-\frac1{2(1-\epsilon)}\langle x\rangle^{1-\epsilon}
}
which satisfies the properties
\eq{
|\grad q_\epsilon|\leq1,\\
\xi^T\grad^2q_\epsilon(x)\xi>\frac\epsilon2\langle x\rangle^{-(1+\epsilon)}|\xi|^2,\\
\frac{1-2\epsilon}{2-2\epsilon}\langle x\rangle\leq q_\epsilon(x)\leq\langle x\rangle.
}
With $\tau$ a time scale to be specified, we define $w\coloneqq q_\epsilon(\tau\omega)$ which obeys the equation
\eqn{
(\dd_t+u\cdot\grad-\Delta)w=\tau\grad q_\epsilon(\tau\omega)\cdot(\omega\cdot\grad u)-\tau^2\tr(\grad\omega^T\grad^2q_\epsilon\grad\omega).
}
Recall that $\omega \coloneqq \mathrm{curl}\,u$ denotes the vorticity vector. Multiplying by a spatial cutoff at length scale $R$ and integrating over $\mathbb R^d$,
\eqn{
\frac {\d}{\d t}\int_\Rd w\psi\leq\int_\Rd (u\cdot\grad\psi+\Delta\psi)w+O(\tau|\grad u|^2)\psi-\frac\epsilon2\tau^2\langle\tau\omega\rangle^{-(1+\epsilon)}|\grad\omega|^2\psi.
}
Let $\tilde\psi$ be an enlarged cutoff function so that $R|\grad\psi|+R^2|\Delta\psi|\leq10\tilde\psi$. We set 
\[
\| f \|_{\LulocR^p} \coloneqq \sup_{B(R)\subset \R^3} \| f \|_{L^p (B(R))}.
\] 
 Integrating in time starting from a $t_0$ to be specified and taking a supremum over the balls,
\eqn{
&\|w\psi(t)\|_{\LulocR^1}+\frac\epsilon2\tau^2\int_{t_0}^t\int_\Rd\langle\tau\omega\rangle^{-(1+\epsilon)}|\grad\omega|^2\psi \,d x\d t\\
&\quad\quad\lesssim \|w(t_0)\|_{\LulocR^1}+\int_{t_0}^t(R^{-2}+R^{-1}\|u\|_{\infty})\|w(t')\|_{\LulocR^1}\d t'+\tau\|\grad u\|_{L_{t,x-\uloc,R}^2}^2.
}
Gr\"onwall's inequality 
\eqn{
\|w(t)\|_{\LulocR^1}\lesssim \left( \|w(t_0)\|_{\LulocR^1}+\tau R A^{O(1)} \right)\exp(R^{-2}|t-t_0|+R^{-1}A^{O(1)}|t-t_0|^\frac12),
}
where $|t-t_0|^{1/2}$ comes from applying the Cauchy-Schwarz inequality in the time integral and by using the energy bound \eqref{localenergy}.
Setting $R=A^{C_1}$ and $\tau=A^{-2C_1}$ for a sufficiently large $C_1$, we find
\eqn{
\|\langle\tau\omega(t)\rangle\|_{\LulocR^1}\lesssim\|\langle\tau\omega(t_0)\rangle\|_{\LulocR^1}.
}
By \eqref{localenergy} and H\"older's inequality, we can find a $t_0\in[1/4,1/2]$ where the right-hand side is bounded by $A^{O(1)}$. Therefore
\eqn{
\int_{t_0}^t\int_\Rd\langle\tau\omega\rangle^{-(1+\epsilon)}|\grad\omega|^2\psi \, \d x\d t&\leq \epsilon^{-1}A^{O(1)}.
}
We use H\"older's inequality with the decomposition
\eqn{
|\grad\omega|^{\frac4{3+\epsilon}}=\big(|\grad\omega|^{\frac4{3+\epsilon}}\langle\tau\omega\rangle^{-2\frac{1+\epsilon}{3+\epsilon}}\big)\langle\tau\omega\rangle^{2\frac{1+\epsilon}{3+\epsilon}}
}
to conclude
\eqn{
\|\grad\omega\|_{\Lut^{4/(3+\epsilon)}([t_0,t]\times\Rd)}\leq\epsilon^{-O(1)}A^{O(1)}.
}
To convert this into a bound on $\grad^2u$, fix a unit ball $B\subset\Rd$ and a cutoff function $\varphi\in C_c^\infty(3B)$ with $\varphi\equiv1$ in $2B$. We decompose $\grad^2u=a+b$ where $a=\grad^2\Delta^{-1}\curl(\varphi\omega)$. Note that $b=\grad f$ where $f=\grad\Delta^{-1}\curl((1-\varphi)\omega)$ is harmonic in $2B$ so for any $p\in[1,\frac43)$,
\eqn{
\|a\|_{L_{t,x}^{p}([t_0,t]\times B)}&\lesssim\|\grad\omega\|_{L_{t,x}^p([t_0,t]\times 3B)}+\|\grad\varphi\|_{L^\infty}\|\omega\|_{\Lut^2([t_0,t]\times\Rd)}\leq\epsilon^{-O(1)}A^{O(1)}
}
and
\eqn{
\|b\|_{L_{t,x}^p([t_0,t]\times B)}&\lesssim\|\grad\Delta^{-1}\curl((1-\varphi)\omega)\|_{L_{t,x}^2([t_0,t]\times2B)}\\
&\lesssim\|\omega^\sharp\|_{L_{t,x}^2([t_0,t]\times\Rd)}+\|\omega^\flat\|_{L_{t,x}^\infty([t_0,t]\times\Rd)}\leq A^{O(1)},
}
where we have used \eqref{localenergy}, H\"older's inequality, \eqref{energy}, and \eqref{flatbounds}.
\end{proof}

\section{Estimates for axisymmetric Navier-Stokes solutions}\label{sec_ests_axisymm_nse}

Here we provide some estimates of classical solutions of \eqref{NSE} that are specific to the axisymmetric assumption on the solutions.

We first note that $u_\theta$ satisfies 
\eqnb\label{utheta_pde}
\left(  \p_t + u \cdot \na - \Delta + \frac1{r^2} \right) u_\theta     +\frac{u_r}{r} u_\theta  =0,
\eqne
which in particular gives that  the swirl $\Theta \coloneqq ru_\theta$ satisfies
 \eq{\label{swirl_pde_repeat}
\Big(\dd_t+\Big(u+\frac2re_r\Big)\cdot\grad-\Delta\Big)\Theta=0
}
in $(\R^3 \setminus \{ r=0 \}) \times (0,T)$. It then follows that, at each time, $(r,x_3)\mapsto u_\theta (r,x_3,t)$ is a continuous function on $\overline{\R_+} \times \R$ with $u_\theta (0,x_3)=0$ for all $x_3$ (see \cite[Lemma~1]{liu_wang} for details). In particular

\eqnb\label{axis_vanishing}
\Theta (0,x_3,t) =0\qquad \text{ for all } x_3\in \R, t\in (0,T).
\eqne

Moreover, since $\omega = \omega_r e_r+ \omega_\theta e_\theta + \omega_3 e_3$ is a smooth vector field we see (also by \cite[Lemma~1]{liu_wang}) that
$\Phi=  \frac{\omega_r}{r}$, $ \Gamma \coloneqq \frac{\omega_\theta}{r}$ (recall~\eqref{phi_gamma_def}) satisfy 
\eqnb\label{phi_gamma_near_origin}
| \Phi (r,x_3 ,t ) |, | \Gamma  (r,x_3 ,t ) | \lec C (x_3,t)
\eqne
for $r\in [0,1]$.

\subsection{H\"older continuity near the axis}\label{sec_holder}

Here we consider the parabolic equation 
\eqnb\label{MV_eq}
\mathcal{M} V \coloneqq \p_t V - \Delta V + b \cdot \nabla V =0
\eqne
in a space-time cylinder 
\[
Q_R (x_0,t_0)\coloneqq B( x_0, R)\times (t_0 - R^2 ,t_0).
\]
We assume that at each point of $Q_R\coloneqq Q_R (0,0)$
\eqnb\label{div_constraint}
\text{either }  \div\, b=0 \quad \text{ or } \quad V=0.
\eqne
We also assume that 
\eqnb\label{N}
 \mathcal N(R )\coloneqq 2+\sup_{R'\leq R  }(R')^{-\alpha}\|b\|_{L_t^{\ell}L_x^q(Q_{R'})}<\infty,
\eqne
where $\alpha\coloneqq \frac 3q+\frac2\ell-1\in[0,1)$. In such setting \cite[Corollary~3.6]{nu_2012} observed that $V$ must be H\"older continuous in the interior of $Q_{R}$, and in the proposition below we state a version of their result in which we quantify the dependence of the H\"older exponent in terms of $\mathcal{N}$. 

\begin{prop}\label{prop_holder}
If $V$ is a Lipschitz solution of \eqref{MV_eq} on $Q_{2R}$ then 
\eqn{
\osc_{B(r)}V(0)\lesssim\left(\frac r R\right)^\gamma\osc_{Q_R}V
}
for all $r\leq R$, where $\gamma=\exp(-{\mathcal N}^{O(1)})$.
\end{prop}
\begin{proof}
See Appendix~\ref{sec_NU}. 
\end{proof}

We note that the swirl $\Theta$ satisfies \eqref{MV_eq} with $b\coloneqq u +2e_r /r $ (recall \eqref{swirl_pde_repeat} above). Moreover $\div \, b = 0$ everywhere except for the axis, since $\div \, u =0$, $\div (e_r /r )=0$ (recall \eqref{div_cylindr}) there. Furthermore, $\Theta=0$ on the axis (recall \eqref{axis_vanishing}), and so the assumption \eqref{div_constraint} holds. Thus Proposition~\ref{prop_holder} shows that $\Theta $ is H\"older continuous in a neighborhood of the axis. We explore this in more detail in the proof of Theorem~\ref{regularity} below, where we quantify $\mathcal{N}$ in terms of the weak-$L^3$ bound $A$ (see Step~1 in Section~\ref{sec_en_ests_cfz} below).

\subsection{Pointwise estimates away from the axis}\label{sec_pointwise}

The following is a more precise version of Proposition 8 in \cite{palasek}.

\begin{prop}[Pointwise bounds away from the axis]\label{pointwisebounds}
Let $u$ solve \eqref{NSE} on $[0,1]$ satisfying \eqref{NSE_axisym} and \eqref{weakL3}. Then for every  $\epsilon \in (0,4/15 )$, we have
\eqn{
|\grad^ju|\leq \left( r^{-1-j}+r^{-\frac13+\epsilon} \right) A^{O_{\epsilon,j}(1)}
}
for each $t\in [1/2,1]$. We also have 
\eqn{
\|u\|_{L^p(\{r\geq1\})}\leq A^{O_p(1)}
}
for each such $t$, and $p\in(3,\infty]$.
\end{prop}

\begin{proof} 
We first pick any $\alpha \in (1/3-\epsilon/2,1/3)$ and $c =c(j)>0$ sufficiently small so that 
\eqnb\label{c_choice}
(1-\alpha + j ) c < \epsilon /2 \quad \text{ and } \quad c< \alpha/ (1-\alpha).
\eqne
We also pick $n=n(j)\in \N$ sufficiently large so that
\eqnb\label{n_choice}
n\geq (2+j )\left( 1 + \frac1c \right) .
\eqne
We set $t_k \coloneqq 1/2 - (1/2)^k$ and we define a sequence of regions $\{x\in\Rd:r\geq R/2\}=\Omega_1\supset\Omega_2\supset\cdots\supset\Omega_n=\{x\in\Rd:r\geq R\}$ such that $\dist(\Omega_i,\Omega_{i+1})\geq R/2n $. \\

Given such a sequence of times we now consider the corresponding Picard iterates $\ulin_k$, $\unlin_k$, for $k\in \{ 0, 1,\ldots , n\}$.\\

\noindent\texttt{Step 1.} We show that
\eq{\label{awayfromtheaxis}
\|P_N \ulin_k (t) \|_{L^\infty (\{ r \geq R /2 \} )},\|P_N \unlin_k (t) \|_{L^\infty (\{ r \geq R /2 \} )}\lesssim R^{-\alpha}N^{1-\alpha}A^{O_k(1)}
}
for all $\alpha\in[0,\frac13)$, $R>0$ and $t\in [t_k,1]$, $k\geq 0$.\\

In fact, we first observe that Lemma~\ref{bernstein} gives that
\eqnb\label{temp00}
\|r^\alpha P_Nu (t) \|_{\infty }\lesssim N^{1-\alpha}\| u(t) \|_{L^{3,\infty}} \lec N^{1-\alpha } A^{O(1)}.
\eqne
Thus, since the first inequality above is valid for any axisymmetric function, it remains to note that the second inequality is also valid for each $\ulin_k$, $\unlin_k$, on $[t_k,1]$, $k\geq 0$. Indeed, the case $k=0$ follows trivially, while the inductive step follows by applying Young's inequality \eqref{young_weak} for weak $L^p$ spaces, and H\"older's inequality \eqref{holder_lorentz} for Lorentz spaces 
\[
\begin{split}
\| \ulin_k (t) \|_{L^{3,\infty }} &\lec \| \Psi (t- t_k ) \|_1 \| u(t_k) \|_{L^{3,\infty }} + \int_{t_k}^t \|\na \Psi (t-t')\|_1 \| (\ulin_{k-1} \otimes  \ulin_{k-1} )(t')\|_{L^{3/2,\infty }}  \d t' \\
&\leq C_k A + C_k \| \ulin_{k-1} \|^2_{L^\infty ([t_{k-1},1 ]; L^{3,\infty })} \int_{t_k}^t (t-t')^{-\frac12 } \d t' \leq A^{O_k(1)}
\end{split}
\]
for $t\in [t_k ,1 ]$, as required, where we also used the heat kernel bounds \eqref{heat_bounds}. \\

\noindent\texttt{Step 2.} We show that the inequality from Step 1 can be improved for $\unlin_k$ for large $k$, namely
\eq{\label{X}
\|P_N\unlin_k\|_{L^\infty([\frac12,1]\times\{r\geq R\})}&\leq NA^{O_k(1)}((RN)^{-(k-1)\alpha}+N^{-(k-1)})
}
for every $k\geq 1$ and $N\in 2^{\N}\cap[100^k\max(1,R^{-1}),\infty)$.\\

We will show that,
\eq{\label{claim}
X_{k,N}\leq N^{-\frac45}A^{O_k(1)}((RN)^{-(k-1)\alpha}+N^{-(k-1)}),
}
for $k\geq 1$ and $N\geq100^k\max(1,R^{-1})$, using induction with respect to $k$, where 
\[
X_{k,N}\coloneqq \|P_N\unlin_k\|_{L^\infty ([t_{k+1},1]; L^{5/3} ( \Omega_k ))}.\]
Then \eqref{X} follows by the local Bernstein inequality \eqref{bernstein_localized}.

As for the base case $k=1$ we note that \eqref{heat} gives that 
\eqn{
\|P_N\unlin_1(t)\|_{{5/3}}&\lesssim\int_{t_1}^t\|P_N\ee^{(t-t')\Delta}\mathbb P\div (u\otimes u)(t')\|_{{5/3}}\d t'\\
&\lesssim\int_{t_1}^t\ee^{-(t-t')N^2/O(1)}N^{\frac65}\|(u\otimes u)(t')\|_{L^{\frac32,\infty}}\d t'\\
&\lesssim N^\frac65\|\ee^{-tN^2/O(1)}\|_{L^1(t_1,1)}\|u\|_{L^{3,\infty}}^2
}
for $t\in[t_1,1]$. Thus
\eq{\label{initialize}
X_{1,N}\leq \|P_N\unlin_1\|_{L^\infty ([t_2,1]; L^{5/3} ) }\leq N^{-\frac45}A^{O(1)},
}
due to H\"older's inequality for Lorentz spaces \eqref{holder_lorentz}.

As for the inductive step, we use the Duhamel formula for $\unlin_k$ (recall \eqref{picard_def}), and the local Bernstein inequality \eqref{bernstein_localized} to obtain 
\[
\begin{split}
\| P_N \unlin_k (t) \|_{L^{5/3} (\Omega_{k})} &\lesssim \int_{t_k}^t \| P_N \ee^{(t-t')\Delta } \mathbb{P} \div (u\otimes u - \ulin_{k-1} \otimes \ulin_{k-1} ) \|_{L^{5/3} (\Omega_k )}   \d t'  \\
&\leq \int_{t_k}^t  N \ee^{-(t-t')N^2/O(1)} \d t' \left(\|  P_N ( u\otimes u - \ulin_{k-1} \otimes \ulin_{k-1} ) \|_{L^\infty ([t_{k},1]; L^{5/3} (\Omega_{k-1}) )} \right.\\
&\hspace{3cm} \left. +(NR)^{-(k-1)\alpha } \|  P_N ( u\otimes u  - \ulin_{k-1} \otimes \ulin_{k-1} )  \|_{L^\infty ([t_k , 1]; L^{5/3})} \right)\\
&\lec N^{-1} \left( \|  P_N ( u\otimes u - \ulin_{k-1} \otimes \ulin_{k-1} ) \|_{L^\infty ([t_{k},1]; L^{5/3} (\Omega_{k-1}) )} + N^{\frac15} (NR)^{-(k-1)\alpha }  A^{O(1)} \right) ,
\end{split}
\]
where we used the weak $L^3$ bound \eqref{weakL3} and Lemma~\ref{bernstein} for the $u\otimes u$ term and \eqref{flatbounds} for the $\ulin_{k-1} \otimes \ulin_{k-1}$ term. Thus we can use the  paraproduct decomposition in the first term on the right-hand side to obtain
\eq{\label{Xdecomposition}
X_{k,N}&\lesssim N^{-1}\|Y_1+\cdots+Y_5\|_{L^\infty ([t_{k},1]; L^{5/3} (\Omega_{k-1}))}+N^{-\frac45}( NR)^{-(k-1)\alpha }A^{O_k(1)},
}
 where
\eqn{
Y_1&\coloneqq 2\sum_{N'\sim N}P_{N'}\unlin_{k-1}\odot P_{\leq N/100}\unlin_{k-1},\\
Y_2&\coloneqq \sum_{N_1\sim N_2\gtrsim N}P_{N_1}\unlin_{k-1}\otimes P_{N_2}\unlin_{k-1},\\
Y_3&\coloneqq \sum_{N_1\sim N_2\gtrsim N}P_{N_1}\ulin_{k-1}\otimes P_{N_2}\unlin_{k-1},\\
Y_4&\coloneqq 2\sum_{N'\sim N}P_{N'}\ulin_{k-1}\odot P_{\leq N/100}\unlin_{k-1},\\
Y_5&\coloneqq 2\sum_{N'\sim N}P_{\leq N/100}\ulin_{k-1}\odot P_{N'}\unlin_{k-1},
}
where we use the notation $a\odot b\coloneqq a\otimes b + b\otimes a$. Using \eqref{awayfromtheaxis},
\eqn{
\|Y_1\|_{L^\infty ([t_k,1];L^{5/3}(\Omega_{k-1}))}&\lesssim\sum_{N'\sim N}X_{k-1,N'}\sum_{N'\lesssim N} R^{-\alpha}(N')^{1-\alpha} A^{O_k(1)}\\
&\lesssim R^{-\alpha}N^{1-\alpha} A^{O_k(1)}\sum_{N'\sim N}X_{k-1,N'}
}
and
\eqn{
\|Y_2\|_{L^\infty ([t_k,1];L^{5/3}(\Omega_{k-1}))}&\lesssim R^{-\alpha}A^{O_k(1)}\sum_{N'\gtrsim  N} (N')^{1-\alpha} X_{k-1,N'}.
}
Moreover, the frequency-localized bounds \eqref{frequencylocalized} for $\ulin_{k-1}$ give that
\eqn{
\|Y_3\|_{L^\infty ([t_k,1]; L^{5/3}(\Omega_{k-1}))}&\lesssim A^{O_k(1)}\sum_{N'\gtrsim N}\ee^{-(N')^2/O_k(1)}N'X_{k-1,N'},
}
and \eqref{sharpbound}, as well as boundedness of $P_{\leq N/100}$ on $L^{5/3}$ give that 
\eqn{
\|Y_4\|_{L^\infty ([t_k,1];L^{5/3}(\Omega_{k-1}))}&\lesssim A^{O_k(1)}\sum_{N'\sim N}\ee^{-(N')^2/O_k(1)}N'\lesssim \ee^{-N^2/O_k(1)}A^{O_k(1)}.
}
Finally, using boundedness of $P_{\leq N/100}$ on $L^\infty$ and \eqref{flatbounds} we obtain
\eqn{
\|Y_5\|_{L^\infty ([t_k,1]; L^{5/3}(\Omega_{k-1}))}&\lesssim A^{O_k(1)}\sum_{N'\sim N}X_{k-1,N'}.
}
Combining these estimates into \eqref{Xdecomposition}, we have shown
\begin{equation}\begin{aligned}\label{induct}
X_{k,N}&\leq A^{O_k(1)}\left( ((RN)^{-\alpha}+N^{-1})\sum_{N'\sim N} X_{k-1,N'} +  N^{-1} R^{-\alpha } \sum_{N' \gtrsim N } (N')^{1-\alpha } X_{k-1,N'}\right.\\
&\quad\left.+N^{-1}\sum_{N'\gtrsim N}\ee^{-(N')^2/O_k(1)}N' X_{k-1,N'}+N^{-\frac45}(NR)^{-(k-1)\alpha }+N^{-1} \ee^{-N^2/O_k(1)}\right).
\end{aligned}\end{equation}

Since the upper bounds on $X_{k-1,N'}$ provided by the inductive assumption \eqref{claim} are comparable for all $N'\sim N$, up to constants depending only on $k$, we thus obtain that
\[\begin{split}
\sum_{N'\sim N} X_{k-1,N'} &\leq A^{O_{k}(1)}N^{-\frac45} \left( (RN)^{-\alpha ({k-2}) } + N^{-{k-2}}  \right) ,\\
R^{-\alpha }\sum_{N' \gtrsim N } (N')^{1-\alpha } X_{k-1,N'} &\leq  A^{O_{k}(1)} R^{-\alpha }\sum_{N' \gtrsim N}(N')^{1-\alpha -\frac45}  \left( (RN')^{-\alpha ({k-2}) } + (N')^{-({k-2})}  \right) \\
&\leq A^{O_{k}(1)} N^{\frac15}   \left( (RN)^{-\alpha ({k-1}) } + N^{-({k-1})}  \right) ,
\end{split} 
\]
where, in the last line we used the fact that $(k-1)(1-\alpha ) - 4/5 <0$ for any $k\geq 2$.  A similar estimate for $\sum_{N'\gtrsim N}\ee^{-(N')^2/O_k(1)}N' X_{k-1,N'}$ now allows us to deduce from \eqref{induct} that 

\[
\begin{split}
X_{k,N} &\leq  N^{-\frac45 }A^{O_k(1)}((RN)^{-(k-1)\alpha}+N^{-(k-1)}),
\end{split}
\]
as required.\\

\noindent\texttt{Step 3.} We prove the claim.\\

We first consider the case $R\geq 100^{n/c}$, and we note that, by \eqref{awayfromtheaxis}
\eqn{
\|P_{N\leq R^{c}}\grad^j\unlin_n\|_{L_{t,x}^\infty([\frac12,1]\times\{r\geq R\})}&\leq\sum_{N\leq R^c}A^{O_n(1)}N^{1-\alpha+j}R^{-\alpha}\leq A^{O_n(1)}R^{-\alpha+(1-\alpha+j)c}\leq A^{O_{n}(1) }R^{-\frac13 +\varepsilon},
}
where we used the choice of $\alpha>1/3-\epsilon/2$ and the first property of our choice \eqref{c_choice} of $c$ in the last inequality.
On the other hand for $N>R^c$ we can use \eqref{X} with $k= n$ to obtain arbitrarily fast decay in $N$. Comparing the terms on the right-hand side of \eqref{X} we see that $N^{-(n-2)}$ dominates $(RN)^{-(n-2)\alpha}$ if and only if $N\leq R^{\alpha /(1-\alpha )}$, which allows us to apply the decomposition
\eqn{
\|P_{N>R^{c}}\grad^j\unlin_n\|_{L_{t,x}^\infty([\frac12,1]\times\{r\geq R\})}&\leq\sum_{R^c<N\leq R^{\alpha/(1-\alpha )}}A^{O_n(1)}N^{-n+2+j}\\
&\quad+\sum_{N>R^{\alpha/(1-\alpha )}}A^{O_n(1)}N^{1+j}(RN)^{-(n-1)\alpha}\\
&\leq A^{O_n (1)} R^{c(-n+2+j)}\\
&\leq A^{O_n(1)}R^{-1-j},
}
where we used the second property of our choice \eqref{c_choice} of $c$ in the second inequality, and the choice \eqref{n_choice} of $n$ in the last inequality.

We now suppose that $R\leq 100^{n/c}$. The low frequencies can be estimated directly from the weak $L^3$ bound \eqref{weakL3},
\eqn{
\|P_{\leq 100^{2n/c}R^{-1}}\grad^ju\|_{L_{t,x}^\infty([\frac12,1]\times\{r\geq R\})}&\lesssim_{n,c} A^{O(1)}R^{-1-j}.
}
On the other hand, for $N> 100^{2n/c}R^{-1}$ we have in particular $N>R^{\alpha/(1-\alpha )}$, which shows that the dominant term on the right-hand side of \eqref{X} is $(RN)^{-(n-2)\alpha}$, and so 
\eqn{
\|P_{>100^{2n/c}R^{-1}}\grad^j\unlin_n (t) \|_{L^\infty (\{r\geq R\})}&\leq\sum_{N>100^{2n/c}R^{-1}}N^{1+j}A^{O_n(1)}(RN)^{-(n-1)\alpha}\leq A^{O_n(1)}R^{-1-j}
}
for every $t\in[1/2,1]$, as desired. As for the estimate for $\ulin$ we use \eqref{flatbounds} to obtain
\[
\| \nabla^j \ulin_n \|_{L^\infty (\{ r\geq R \} )} \leq R^{-1/3+\epsilon } \|r^{1/3-\epsilon} \nabla^j \ulin_n \|_{\infty } \lec_\epsilon  R^{-1/3+\epsilon } A^{O_{\epsilon , j} (1)},
\]
as needed.

The estimate for $\| u \|_{L^p (\{ r \geq 1\} )}$ follows by an $L^p$ analogue of Step 1, as well as applying the $X_{k,N}$ estimates \eqref{claim} in the $L^p$ variant of Step 3. 
\end{proof}

\subsection{A Poisson-type estimate on $u_r/r$}\label{sec_poisson}

Here we discuss how derivatives of $u_r/r$  can be controlled by $\Gamma$ using the representation \eqref{ur_trick},
\eqnb\label{ur_trick1}
\frac{u_r}r = \Delta^{-1} \p_3 \Gamma - 2 \frac{\p_r }{r} \Delta^{-2} \p_3 \Gamma,
\eqne
see \cite[p. 1929]{cfz_2017}, which will be an essential part of our $x_3$-$\uloc$ energy estimates for $\Phi$ and $\Gamma $ (see Proposition~\ref{prop_energy} below).

\begin{lem}[The $L^2_{3-\uloc}$ estimate on $u_r/r$]\label{lem_localgamma}
\eqnb\label{localgamma_claim_repeat}
\Big\|\nabla \p_r \frac{u_r}r\Big\|_{\Luloc^2}+ \Big\|\nabla \p_3 \frac{u_r}r\Big\|_{\Luloc^2}\lesssim\|\Gamma\|_{\Luloc^2}+\|\grad \Gamma\|_{\Luloc^2}.
\eqne
\end{lem}

 A version of the above estimate without the localization in $x_3$ has appeared in \cite[Lemma~2.3]{cfz_2017}. As mentioned in the introduction, the localization makes the estimate much more challenging, particularly due to the bilaplacian term in \eqref{ur_trick1}. \\
 
  In order to prove Lemma~\ref{lem_localgamma} we note that, since 
  \[
\frac{\p_r }r = \Delta' - \p_{rr},
\]
\eqref{ur_trick1} gives that
\eqnb\label{ur_r}
\frac{u_r}r = -\Delta^{-1} \p_3 \Gamma + 2(\p_{rr} -\Delta'  )  \Delta^{-2} \p_3 \Gamma.
\eqne
Thus, since $| \na \p_3 \frac{u_r}{r} | = |(\p_r \p_3 \frac{u_r}{r} , \p_3 \p_3 \frac{u_r}{r})| $ (and similarly for $|\na \p_r \frac{u_r}r |$), we can use \eqref{drr_main} and \eqref{d4_in_r_x_3} to observe that
\[
\begin{split}
\left|\na \p_3 \frac{u_r}{r} \right|+ \left|\na \p_r \frac{u_r}{r} \right| &\lec |D^2_{r,x_3} \Delta^{-1} \p_3 \Gamma | + |D^2_{r,x_3} (\p_{rr} - \Delta')\Delta^{-2} \p_3 \Gamma |\\
&\lec |\nabla \Gamma | + | D^2 \Delta^{-1} \nabla' \Gamma | + | D^4  \Delta^{-2} \nabla' \Gamma |,
\end{split}
\] 
where we used $\p_{33} = \Delta - \Delta'$ in the last line. In particular, each of the terms on the right-hand side involves at least one derivative in the horizontal variables.
Thus, in order to estimate the left-hand side of \eqref{localgamma_claim_repeat} it suffices to find suitable bounds on the last two terms, which we achieve in Lemmas~\ref{lem_L}--\ref{lem_doubleL} below. Their claims give us \eqref{localgamma_claim_repeat}, as required.\\

\begin{lem}\label{lem_L}
If $f= \Delta^{-1} \nabla' \Gamma $ then
\[
\| D^2 f \|_{L^2_{3-\uloc}} \leq \| \Gamma \|_{L^2_{3-\uloc}} +  \| \nabla' \Gamma \|_{L^2_{3-\uloc}} .
\]
\end{lem}
\begin{proof}
Let $I(x)$ denote the kernel matrix of $D^2 (-\Delta )^{-1}$. We have that 
\[
|\nabla^j I(x) | \leq \frac{C}{|x|^{3+j}} \qquad \text{ for } j=0,1,
\]
and
\[\begin{split}
D^2 f(x)& = \mathrm{p.v.}\int_{\R^3} I( x-y ) \nabla' \Gamma (y)   \d y \\
&= \mathrm{p.v.} \int_{\R^3} \nabla' \Gamma (y) \tilde \phi (y_3) I(x-y ) \d y +   \mathrm{p.v.} \int_{\R^3}  \Gamma (y) (1-\tilde \phi (y_3) ) \nabla' I( x-y ) \d y  \\
& =: f_1 (x) + f_2 (x).
\end{split}
\]
 The Calder\'on-Zygmund inequality (see \cite[Theorem~B.5]{nse_book}, for example) gives that 
\[
\| f_1 \|_{L^2_{3-\uloc}}  \leq \| f_1 \|_{L^2}  \lec  \| \nabla' \Gamma \, \tilde \phi \|_{L^2} \lec  \| \nabla' \Gamma  \|_{L^2_{3-\uloc}}.
\]
Moreover, noting that $\int_{\R^2 } \frac{ \d x_1 \, \d x_2}{(a^2 +x_1^2 +x_2^2)^2} = C a^{-2} $, we can use Young's inequality for convolutions to obtain  
\[
\begin{split}
\| f_2 (\cdot , x_3 ) \|_{L^2} &\leq  \int_{\R} \frac{ \| \Gamma (\cdot , y_3 ) \|_{L^2} (1-\tilde \phi (y_3) )}{|x_3-y_3|^2} \d y_3 \\
&\leq \sum_{j\geq 1} \int_{\{ | x_3 - y_3 | \in (j,j+1) } \frac{ \| \Gamma (\cdot , y_3 ) \|_{L^2} (1-\tilde \phi (y_3) )}{|x_3-y_3|^2} \d y_3 \\
& \leq \sum_{j\geq 1} j^{-2} \int_{\{ | x_3 - y_3 | \in (j,j+1) } \| \Gamma (\cdot , y_3 ) \|_{L^2}  \d y_3 \\
& \leq \| \Gamma \|_{L^2_{3-\uloc}}.
\end{split} 
\]
integration in $x_3$ over $\supp\, \phi$ finishes the proof.
\end{proof}
For the bilaplacian term in \eqref{ur_r} one needs to work harder:
\begin{lem}\label{lem_doubleL}
Let $f=D^4 \Delta^{-2} \nabla'  \Gamma $. Then
\[
\| f \|_{L^2_{3-\uloc}} \leq \| \Gamma \|_{L^2_{3-\uloc}} +  \| \nabla \Gamma \|_{L^2_{3-\uloc}} .
\]
\end{lem}
\begin{proof} We have that
\[
f(x) = \mathrm{p.v.}\int_{\R^3} \mathrm{p.v.}\int_{\R^3} \nabla' \Gamma (z)  I(x-y) I(y-z)   \d z \, \d y.
\]
Recalling that $\tilde \phi = \sum_{|j|\leq 10 } \phi_j$, and $\tilde {\tilde \phi } =  \sum_{|j|\leq 20 } \phi_j$ we use the partition of unity,
\[
\begin{split} 
1&= \tilde  {\tilde \phi } (z_3) +  (1- \tilde {\tilde \phi } (z_3 ) ) {\tilde \phi } (y_3) +  \sum_{\substack { |j|>10 \\ |k|>20 } } \phi_j (y_3) \phi_k (z_3) \\
&=  \tilde  {\tilde \phi } (z_3) +  (1- \tilde {\tilde \phi } (z_3 ) ) {\tilde \phi } (y_3)\\
&+  \sum_{|j| > 10 } \phi_j (y_3) \left( \sum_{\substack{|k|>20\\|k-j |\leq 10}}\phi_k (z_3) + \sum_{\substack{|k|>20\\ |k-j|>10 \\ k\leq j/2}  }\phi_k (z_3)+ \sum_{\substack{|k|>20\\ |k-j|>10 \\ j/2<k\leq 2j}  }\phi_k (z_3)+ \sum_{\substack{|k|>20\\|k-j|>10 \\ k> 2j}  }\phi_k (z_3) \right) ,
\end{split}
\] 
to decompose $ f$ accordingly,
\[
\begin{split}
f(x) &=  \mathrm{p.v.} \int_{\R^3} \mathrm{p.v.} \int_{\R^3}  \nabla' \Gamma (z) \tilde {\tilde \phi} (z_3) I( x-y) I(y-z)  \d y \, \d z \\
&\quad +\mathrm{p.v.}\int_{\R^3}I(x-y)\tilde \phi (y_3)  \mathrm{p.v.} \int_{\R^3}  \nabla' \Gamma (z)  (1-\tilde {\tilde \phi} (z_3))  I(y-z)  \d z \, \d y \\
&\quad + \mathrm{p.v.} \int_{\R^3}I(x-y)\sum_{|j|>10 } \phi_j (y_3)  \mathrm{p.v.} \int_{\R^3}  \nabla' \Gamma (z) \sum_{\substack{|k|>20 \\ |k-j | \leq 10 }} \phi_k (z_3 )   I(y-z)  \d z \, \d y \\
&\quad + \mathrm{p.v.} \int_{\R^3}I(x-y)\sum_{|j|>10 } \phi_j (y_3)  \mathrm{p.v.} \int_{\R^3}  \nabla' \Gamma (z) \sum_{\substack{|k|>20 \\ |k-j | > 10\\ k\leq j/2 }} \phi_k (z_3 )   I(y-z)  \d z \, \d y \\
&\quad +  \mathrm{p.v.}\int_{\R^3}I(x-y)\sum_{|j|>10 } \phi_j (y_3)  \mathrm{p.v.} \int_{\R^3}  \nabla' \Gamma (z) \sum_{\substack{|k|>20 \\ |k-j | > 10 \\ j/2 < k \leq 2j  }} \phi_k (z_3 )   I(y-z)  \d z \, \d y \\
&\quad +  \mathrm{p.v.}\int_{\R^3}I(x-y)\sum_{|j|>10 } \phi_j (y_3)  \mathrm{p.v.} \int_{\R^3}  \nabla' \Gamma (z) \sum_{\substack{|k|>20 \\ |k-j | > 10 \\ k >2j  }} \phi_k (z_3 )   I(y-z)  \d z \, \d y \\
&=:  f_1 (x) + f_2 (x) + f_3 (x) + f_4 (x) + f_5 (x) +  f_6 (x) .
\end{split} 
\]
Clearly $ f_1$ involves localization of $ \nabla'  \Gamma$ in $z_3$, and so we can use the Calder\'on-Zygmund inequality twice to obtain
\[
 \| f_1 \|_{L^2}  \lec \| \nabla  \Gamma \|_{L^2_{3-\uloc}}.
\]
As for $f_2$ we integrate by parts in the $z$-integral (note that this does not conflict with the principal value, as the singularity has been cut off, and the far field has sufficient decay) and apply the Calder\'on-Zygmund estimate in $x$ to obtain
\[
\begin{split}
 \| f_2 \|_{L^2}  &\lec \left\| \tilde \phi (y_3 ) \int_{\R^3} \frac{ |\Gamma (z)| (1-\tilde {\tilde \phi } (z_3))}{ | y-z |^4 } \d z   \right\|_{L^2} \lec  \sup_{y_3 \in \supp \, \tilde \phi } \left\| \int_{\R^3} \frac{ |\Gamma (z)| (1-\tilde {\tilde \phi } (z_3))}{ | y-z |^4 } \d z   \right\|_{L^2_{y'}}\\
 & \lec  \sup_{y_3 \in \supp \, \tilde \phi }  \int_{\R} \frac{ \|\Gamma (\cdot , z_3 ) \|_{L^2} (1-\tilde {\tilde \phi } (z_3))}{ | y_3-z_3 |^2 } \d z_3  \\
 & \lec  \sup_{y_3 \in \supp \, \tilde \phi }   \sum_{j \geq 1} j^{-2} \int_{|z_3 - y_3 |\in (j,j+1)} \| \Gamma (\cdot , z_3 ) \|_{L^2_{z'}} \d z_3 \lec \| \Gamma \|_{L^2_{3-\uloc}},
 \end{split}
\]
where we used Young's inequality in the second line (as in the lemma above).

As for $f_3$, we integrate by parts in $z$ and then in $y$ to obtain
\[
|f_3 (x) | \lec \sum_{|j|>10}  \int_{\R^3} \frac{ \phi_j (y_3) }{|x-y|^4} \left| \mathrm{p.v.} \int_{\R^3}   \Gamma (z) \sum_{\substack{|k|>20 \\ |k-j | \leq 10 }} \phi_k (z_3 )   I(y-z)  \d z \right| \d y .
\]
We note that the integration by parts is justified as
\[
f_3 = D^2 (-\Delta)^{-1} \left((1-\sum_{|j|\leq 10 } \phi_j (y_3) ) D^2 (-\Delta )^{-1} \left( \nabla' \Gamma (1-\sum_{k\in I} \phi_k (z_3) \right)  \right),
\]
where $I\coloneqq \{ -20 ,\ldots , 20 \} \cup \{ j-10 , \ldots , j+10\} $ is a finite index set. Thus, the operation of integration by parts above is equivalent to moving  $\nabla'$  outside of the outer brackets, which in turn holds since the sums do not depend on $x'$ and $\nabla'$ commutes with other differential symbols.

Thus, using Young's inequality in $x'$
\[
\begin{split}
 \| f_3 (\cdot , x_3) \|_{L^2_{x'}}  & \lec \sum_{|j|>10}  \int_{\R} \frac{ \phi_j (y_3) }{|x_3-y_3|^2} \left\| \mathrm{p.v.} \int_{\R^3}   \Gamma (z) \sum_{\substack{|k|>6 \\ |k-j | \leq 2 }} \phi_k (z_3 )   I(y-z)  \d z \right\|_{L^2_{y'}} \d y_3 \\
  & \lec  \sum_{|j|>2} j^{-2}  \left\| \mathrm{p.v.} \int_{\R^3}   \Gamma (z) \sum_{\substack{|k|>20 \\ |k-j | \leq 10 }} \phi_k (z_3 )   I(y-z)  \d z \right\|_{L^2_{y}} \\
  & \lec  \sum_{|j|>10} j^{-2}  \left\|    \Gamma (z) \sum_{\substack{|k|>20 \\ |k-j | \leq 10 }} \phi_k (z_3 )  \right\|_{L^2}  \lec \| \Gamma \|_{L^2_{3-\uloc}}
 \end{split}
\]
for each $x_3 \in \supp\, \phi$, where we applied the Cauchy-Schwarz inequality (in $y_3$) in the second line.

As for $f_4$ we note that 
\[
|y_3-z_3 | \geq |y_3| - |z_3| \geq (j-1) - (k+1) \geq \frac{j}{2} -2 \geq (j+2)/4 \geq  (|y_3 | +1 ) /4 \geq |y_3 - x_3|/4.
\]
Thus, we can integrate by parts in $z$ to obtain
\[
|f_4 (x) | \leq \int_{\R^3} \int_{\R^3 \cap \{  |y_3 -z_3 | \geq |x_3 - y_3 |/4 \} } \frac{|\Gamma (z) | (1-\tilde{\phi } (y_3 )) (1-\tilde{\phi }(y_3-z_3 ))}{|x-y|^3 |y-z|^4 } \d z \, \d y.
\]
Hence, applying Young's inequality in $x'$ and then in $y'$ we obtain
\eqnb\label{double_young_ex}
\begin{split}
\| f_{4} (\cdot , x_3) \|_{L^2} &\leq \int_{\R } \left\| \int_{\R^3\cap \{  |y_3 -z_3 | \geq |x_3 - y_3 |/4 \}}  \frac{ \Gamma (z) (1-\tilde \phi (y_3) )  (1-\tilde \phi (y_3 - z_3 ))  }{|y-z|^4}  \d z \right\|_{L^2_{y'}} \\
&\hspace{8cm} \cdot \underbrace{ \int_{\R^2} \frac{\d x_1 \, \d x_2}{\left( |x_3 - y_3 |^2 + x_1^2 +x_2^2 \right)^{3/2} }  }_{=C|x_3-y_3 |^{-1}}\d y_3\\
&\lec \int_{\R } \int_{\R \cap \{  |y_3 -z_3 | \geq |x_3 - y_3 |/4 \}} \frac{\| \Gamma (\cdot , z_3 ) \|_{L^2} (1-\tilde \phi (y_3) )  (1-\tilde \phi (y_3 - z_3 ))  }{|x_3 - y_3 | \, |y_3 - z_3 |^2} \d z_3 \, \d y_3 .
\end{split}
\eqne
Hence 
\[\begin{split}
\| f_4 (\cdot , x_3 )\|_{L^2}& \leq \int_{\R} \frac{1-\tilde \phi (y_3 ) }{|x_3 - y_3|^{3/2}}  \left( \sum_{j\geq 1 } \int_{\{|y_3 - z_3 | \in (j,j+1) \} }  \frac{ \| \Gamma (\cdot , z_3 ) \|_{L^2}}{  | y_3 - z_3 |^{3/2}} \d z_3 \, \right)  \d y_3 \\
& \lec \| \Gamma \|_{L^2_{3-\uloc}} \int_{\R} \frac{1-\tilde \phi (y_3 ) }{|x_3 - y_3|^{3/2}} \d y_3 \lec \| \Gamma \|_{L^2_{3-\uloc}}.
\end{split}
\]

As for $f_5$ we have
\[
\frac14 \leq \frac{|x_3 -y_3 |}{|x_3 - z_3 | } \leq 4 ,
\]
since 
\[
|x_3 - y_3 | \leq |y_3 | + |x_3 | \leq j+2 \leq 2j-8 \leq 4k - 8 \leq 4 (|z_3 | - |x_3 | )\leq 4 | x_3 - z_3 |
\]
and 
\[
|x_3 - z_3 | \leq |z_3| + |x_3| \leq k+2 \leq 2j+2 \leq 4(j-2) \leq 4 (| y_3 | - |x_3 | ) \leq 4 |x_3 - y_3 |.
\]
In particular, the triangle inequality gives that
\[
|y_3 - z_3  | \leq 5|x_3 -z_3 |.
\]
Thus we can integrate by parts twice (in $z$ and then in $y$, so that the derivative falls on $I(x-y)$), and then use Young's inequality twice  (as in \eqref{double_young_ex} above) and Tonelli's Theorem to obtain
\[
\begin{split}
\| f_5 (\cdot , x_3 ) \|_{L^2}&\leq  \int_\R \int_{\{ |x_3 - y_3 | /4\leq |x_3-z_3 | \leq 4  | x_3 - y_3 | \} } \frac{ \| \Gamma (\cdot , z_3 ) \|_{L^2} (1- \tilde \phi (y_3 - z_3 )) (1-\tilde {\tilde \phi} (z_3 )) }{|x_3-y_3 |^2 | y_3 - z_3 |} \d z_3 \, \d y_3 \\
&\leq \int_\R \frac{\| \Gamma (\cdot , z_3 ) \|_{L^2} (1-\tilde {\tilde \phi} (z_3 ))}{|x_3 - z_3 |^2 }  \int_{\{ |y_3-z_3 | \leq 5 |x_3 - x_3 | \} } \frac{ 1- \tilde \phi (y_3 - z_3 )  }{ | y_3 - z_3 |} \d y_3 \, \d z_3 \\
&\lec \int_\R \frac{\| \Gamma (\cdot , z_3 ) \|_{L^2} (1-\tilde {\tilde \phi} (z_3 ))}{|x_3 - z_3 |^2 } \log (5|x_3-z_3|)\d z_3 \\
&\lec \sum_{j\geq 1 } \int_{|z_3-x_3| \in (j,j+1)} \frac{\| \Gamma (\cdot , z_3 ) \|_{L^2} }{|x_3-z_3|^2 } \log (5|x_3-z_3|) \d z_3 \\
 &\lec \sum_{j\geq 1 } j^{-2}\log (5j) \| \Gamma  \|_{L^2_{3-\uloc}} \lec \| \Gamma  \|_{L^2_{3-\uloc}}.
 \end{split} 
\]
Finally, for $f_6$ we observe that 
\[
\frac14 \leq \frac{ |x_3 - z_3 |}{|y_3 - z_3 |} \leq 4,
\]
since
\[
|y_3-z_3| \geq |z_3| - |y_3|  \geq k - j -2 > \frac{k-8 }2\geq \frac{k+2}{4} \geq \frac{|x_3| + |z_3| }4 \geq  \frac{|x_3- z_3|}4 
\]
and
\[
|y_3-z_3 | \leq | y_3| + |z_3| \leq j+k +2 \leq \frac{3k+4}{2} \leq 4(k-2) \leq 4(|z_3| - |x_3| )\leq 4|x_3-z_3|.
\]
In particular, the triangle inequality gives that 
\[
| x_3 - y_3 | \leq 5 |x_3 -z_3|.
\]
Thus, similarly to the  case of $f_5$ (although without integrating by parts in $y$), we apply Young's inequality twice, and Tonelli's Theorem to obtain
\[
\begin{split}
\| f_6 (\cdot , x_3 ) \|_{L^2} &\leq  \int_\R \int_{ \R } \frac{ \| \Gamma (\cdot , z_3 ) \|_{L^2} (1-\tilde \phi (y_3)) (1-\tilde {\tilde \phi} (z_3 )) }{|x_3-y_3 | | y_3 - z_3 |^2} \d z_3 \, \d y_3 \\
 &\leq \int_\R \frac{\| \Gamma (\cdot , z_3 ) \|_{L^2} (1-\tilde {\tilde \phi} (z_3 )) }{|x_3-z_3|^2 } \int_{\{ \frac14 |z_3 - x_3 |\leq |y_3 - z_3 | \leq 4 |z_3 - x_3 | \} } \frac{ 1-\tilde \phi (y_3)  }{|x_3-y_3 | } \d y_3 \, \d z_3 \\
 &\leq \int_\R \frac{\| \Gamma (\cdot , z_3 ) \|_{L^2} (1-\tilde {\tilde \phi} (z_3 )) }{|x_3-z_3|^2 } \int_{\{ 1\leq |x_3-y_3 | \leq 5 |x_3-z_3 | \}  } \frac{ 1  }{|x_3-y_3 | } \d y_3 \, \d z_3 \\
 &\lec  \int_\R \frac{\| \Gamma (\cdot , z_3 ) \|_{L^2} (1-\tilde {\tilde \phi} (z_3 )) }{|x_3-z_3|^2 } \log(5|x_3-z_3|) \d z_3  \\
 &\lec \sum_{j\geq 1 } \int_{|z_3-x_3| \in (j,j+1)} \frac{\| \Gamma (\cdot , z_3 ) \|_{L^2}\log(5|x_3-z_3|) }{|x_3-z_3|^2 }  \d z_3 \\
 &\lec \sum_{j\geq 1 } \log (5j )j^{-2}\| \Gamma  \|_{L^2_{3-\uloc}} \lec \| \Gamma  \|_{L^2_{3-\uloc}}
 \end{split} 
\]
for $x_3 \in \supp\, \phi$. Integration of the squares of the above estimates for $f_3,f_4,f_5,f_6$ gives the claim.
\end{proof}

\section{Energy estimates for $\omega/r$}\label{sec_en_ests_cfz}

In this section, we assume the weak $L^3$ bound \eqref{weakL3} on the time interval $[0,1]$ and prove an energy bound for $\Phi^2 + \Gamma^2$ at time $1$, that is we prove the following.

\begin{prop}[An $L^2_{3-\uloc}$ energy estimate for $\Phi$ and $\Gamma$]\label{prop_energy}
Let $u$ be a classical solution of \eqref{NSE} satisfying the weak $L^3$ bound \eqref{weakL3} on $[0,1]$. Then
\eqnb\label{claim_of_prop}
\|\Phi (1)\|_{\Luloc^2(\Rd)} + \|\Gamma(1)\|_{\Luloc^2(\Rd)} \leq \exp\exp A^{O(1)}.
\eqne
\end{prop}
Recall \eqref{l2_3uloc_def} that $\|\cdot \|_{L^p_{3-\uloc} }\coloneqq \sup_{z\in\mathbb R}\|\cdot \|_{L^p(\mathbb R^2\times[z-1,z+1])}$.  We note that we will only use (in \eqref{B} below) the bound on $\Gamma$. 
\begin{proof}
We fix a cutoff function $\phi\in C_c^\infty((-1,1);[0,1])$ such that $\phi\equiv1$ in $[-1/2,1/2]$, and we define the translate 
\[
\phi_z(y)\coloneqq \phi(y-z).
\]
Clearly, we have the pointwise inequality
\eqn{
\phi_z',\phi_z''\lesssim\sum_{i=-2}^2\phi_{z+i}.
}
We will consider the energies
\eqn{
&E(t)\coloneqq \sup_{z\in\mathbb R}E_z(t),\quad E_z(t)\coloneqq \frac12\int_\Rd(\Phi(t,x)^2+\Gamma(t,x)^2)\phi_z(x_3)\d x,\\
&F(t)\coloneqq \sup_{z\in\mathbb R}F_z(t),\quad F_z(t)\coloneqq \int_{t_0}^t\int_\Rd(\grad\Phi(s,x)^2+\grad\Gamma(s,x)^2)\phi_z(x_3)\d x\, \d s
}
for $t\in [t_0,1]$, where $t_0\in[0,1]$ will be chosen in Step 3 below. Given $z\in\mathbb R$, we multiply the equations \eqref{phi_gamma_pdes} by $\phi_z \Gamma$ and $\phi_z \Phi$, respectively, and integrate to obtain, at a given time $t$,
\eqnb\label{en_ineq}
\begin{split}
E_z'&\leq\int_\Rd\Big(-(|\grad\Phi|^2+|\grad\Gamma|^2)\phi_z+\frac12(\Phi^2+\Gamma^2)(u_z\phi_z'+\phi_z'')
\\&\quad\quad\quad+(\omega_r\dd_r+\omega_3\dd_3)\frac{u_r}r\Phi\phi_z-2r^{-1}u_\theta\Phi\Gamma\phi_z\Big)\d x\\
&=:-F_z'(t) + I_1+I_2+I_3.
\end{split}
\eqne
The second term on the right hand side can be bounded directly,
\eqnb\label{i1}
I_1\lesssim(1+\|u_z\|_{L_x^\infty(\Rd)})E(t).
\eqne
The remaining terms $I_2,I_3$ are more challenging. In order to estimate them, as well as choose $t_0$ and deduce the claim \eqref{claim_of_prop}, we follow the steps below.  \\

\noindent\texttt{Step 1.} We use the H\"older estimate (Proposition~\ref{prop_holder}) to show that $|\Theta|\leq r^\gamma A^{O(1)}$ whenever $r\leq\frac12$ and $t\in [3/4,1]$, where $\gamma=\exp(-A^{O(1)})$.\\

To this end we note that, due to incompressibility,  $\div(u+\frac2re_r)=4\pi\delta_{\{x'=0\}}$, which enables us to apply Proposition~\ref{prop_holder} to the equation for the swirl $\Theta$ (recall \eqref{swirl_pde}).

 Moreover, in the notation of Proposition~\ref{prop_holder}, for every $R<\frac12$, $t_0\in[\frac12,1]$ and $x_0\in(0,0)\times\mathbb R$ (i.e., on the $x_3$-axis),
\eqn{
R^{-\frac45}\|u+\frac{e_r}r\|_{L_t^\infty L_x^\frac53(Q((t_0,x_0),R))}&\lesssim R^{-\frac12}\|u\|_{L_t^\infty \Lu^2([t_0-R^2,t_0]\times\Rd)}+1\leq A^{O(1)},
}
by H\"older's inequality and \eqref{localenergy} applied on the timescale $R^2$. (In particular note that each scale $R$ leads to a different decomposition $u=\ulin_n+\unlin_n$, but they all obey the same bounds up to being suitably rescaled.) Thus, for every $r\in(0,1/2)$,
$\osc_{B(x_0,r)}\Theta(t_0)\lesssim r^\gamma\osc_{Q(1/2)}\Theta$ for $r\in (0,1/2)$, which implies the claim.\\

 \noindent\texttt{Step 2.} We show that
 \[
 \int_{t_0}^t \left| I_2 + I_3 \right| \lec \frac12 F(t) + r_0^{-10} + \int_{t_0}^t G E 
 \]
 for each $t_0\in [3/4,1]$ and $t\in [t_0,1]$, where
\eqnb\label{r_0_def}
r_0 \coloneqq \ee^{-\gamma^{-2}},
\eqne
$\gamma = \exp ( -A^{O(1)} )$ is given by Step 1, and
\[
G \coloneqq r_0^{-3} + \| u\|_\infty + \| D^2 u \|_{L^{5/4}_{\uloc}} + \| \nabla u \|_{L^2_{\uloc}}
\]
at each $t'\in [t_0,t]$.\\

To this end, we proceed similarly to \cite{cfz_2017}. Using integration by parts, we compute
\eqn{
I_2&=2\pi\int_\mathbb R\int_0^\infty(-\dd_3u_\theta\dd_r\frac{u_r}r\Phi+\frac{\dd_r(ru_\theta)}r\dd_3\frac{u_r}r\Phi)\phi_z(x_3)r\, \d r \,\d x_3\\
&=\int_\Rd u_\theta(\dd_r\frac{u_r}r\dd_3\Phi\phi_z-\dd_3\frac{u_r}r\dd_r\Phi\phi_z+\dd_r\frac{u_r}r\Phi\phi_z')\\
&=:I_{2,1}+I_{2,2}+I_{2,3}.
}
Let us further decompose $I_{2,i}=I_{2,i,\iin}+I_{2,i,\out}$ ($i=1,2,3$) by writing 
\[
\int = \int_{\{ r < r_0 \} } + \int_{\{ r\geq r_0 \} }.
\]
We decompose
\eqn{
I_{2,1,\iin}&=I_{2,1,\iin,1}+I_{2,1,\iin,2},
}
where
\eqn{
I_{2,1,\iin,1}\coloneqq  \int_{ \{  r<r_0 \} }  u_\theta\left(\fint_\Omega\dd_r\frac{u_r}r\right)\dd_3\Phi\phi_z 
}
and $\Omega\coloneqq \{ x' \colon r <1 \} \times\supp\phi_z$. We compute using H\"older's inequality and Sobolev embedding
\eqn{
\bigg|\int_\Omega\dd_r\frac{u_r}r\bigg|&\leq\|r^{-1}\dd_ru_r\|_{L^1(\Omega)}+\|r^{-2}u_r\|_{L^1(\Omega)}\\
&\lec \|r^{-1}\|_{L^{15/8} (\Omega )}\|\grad u\|_{L^{15/7}(\Omega)}\lesssim\|\grad^2u\|_{L^{5/4}(\Omega)}+\|\grad u\|_{L^2(\Omega)} \lec G.
}
Thus, integrating by parts, and applying H\"older's inequality in Lorentz spaces \eqref{holder_lorentz}, and Young's inequality, we obtain
\eqn{
|I_{2,1,\iin,1}|&\leq G \int_{B(r_0)\times\mathbb R} \left( r\Phi^2\phi_z+|u_\theta\Phi\phi_z'| \right) \d x\\
&\lesssim G (r_0E+\|u_\theta\|_{L_x^{3,\infty}(\Rd)}\|\Phi\|_{L_x^2(\Omega)}|\Omega|^\frac16)\\
&\lesssim G(E+A^{O(1)}).
}
As for $I_{2,1,\iin,2}$ we note that $p=2(1-\gamma)/(1-2\gamma)$ is such that $p-2=2\gamma /(1-2\gamma ) \geq \gamma$ and so we can use the quantified Hardy inequality (Lemma~\ref{lem_hardy}) to obtain, for  we estimate for $t\in[\frac12,1]$,
\eqn{
|I_{2,1,\iin,2}|&\lesssim\|r^{\frac3p-\frac12}u_\theta\|_{ L^{\left( \frac12-\frac1p\right)^{-1}}(\{r\leq r_0\}\cap  \supp\phi_z)}\left\| r^{-\frac3p+\frac12}(\dd_r\frac{u_r}r-\fint_\Omega\dd_r\frac{u_r}r)\phi_z^\frac12\right\|_{L^p(\{r\leq1\})}\|\dd_3\Phi\phi_z^\frac12\|_{2}\\
&\lesssim \gamma^{-O(1)}r_0^{\gamma/3}\left\|\grad\dd_r\frac{u_r}r\right\|_{\Luloc^2 } \|\grad\Phi\phi_z^\frac12\|_{2}\\
&\leq \ee^{-\gamma^{-1}/4}(\|\grad\Gamma\|_{\Luloc^2(\Rd)}+\|\Gamma\|_{\Luloc^2(\Rd)})\|\grad\Phi\phi_z^\frac12\|_{2},
}
where we have also applied  Poincar\'e's inequality and our choice \eqref{r_0_def} of $r_0$. Thus
\[ \int_{t_0}^tI_{2,1,\iin,2}\leq\frac1{20}F(t)+\int_{t_0}^tE.
\]
An analogous argument, in which ``$\p_r$'' and ``$\p_3$'' are switch, gives us the same bound for $I_{2,2,\iin,2}$. As for $I_{2,2,\iin,1}$, we integrate by parts, and apply H\"older's inequality for Lorentz spaces \eqref{holder_lorentz}, and Young's inequality, to obtain
\eqn{
|I_{2,2,\iin,1}|&\leq\bigg|\fint_\Omega\dd_3\frac{u_r}r\phi_z\bigg|\int_{\{r\leq r_0\}\cap\supp\phi_z} |u_\theta\dd_r\Phi|\\
&\lesssim\bigg|\fint_\Omega\frac{u_r}r\phi_z'\bigg|\|u_\theta\|_{L^{3,\infty}}\|\grad\Phi\|_{L_x^2(\supp\phi_z)}r_0^\frac13\\
&\lesssim\sum_{i=-2}^2\|\grad u\|_{L^1(\Omega)}A(F_{z+i}')^\frac12r_0^\frac13\\
&\lesssim G A r_0^{1/3} \left( \sum_{i=-2}^2F_{z+i}'\right)^\frac12,
}
which, thanks to the smallness of $r_0 = \exp(-\exp (A^{O(1)}))$ (recall \eqref{r_0_def}), gives that 
\eqn{
\int_{t_0}^t|I_{2,2,\iin,1}|\leq\frac1{20}F(t)+(t-t_0).
}

We similarly decompose $I_{2,3,\iin}=I_{2,3,\iin,1}+I_{2,3,\iin,2}$ to find
\eqn{
|I_{2,3,\iin,1}|=\left| \fint_\Omega\dd_r\frac{u_r}r \right| \,\, \left| \int_{\{r\leq r_0\}}u_\theta\Phi\phi_z' \right|&\lesssim(\|\grad u\|_{L^{2}(\Omega)}+\|\grad^2u\|_{L^{5/4}(\Omega)})AE^\frac12r_0^\frac13 \\&\lec  G (E+1),
}
where we have used Lemma~\ref{lem_hardy} and change of variables, the pointwise estimate $|u_r/r|\leq|\grad u|$, and H\"older's inequality to bound
\eqn{
\left|\fint_\Omega\dd_r\frac{u_r}r\right|&\lesssim\int_{z-10}^{z+10}\int_0^1\left(|\dd_ru_r|+\frac{|u_r|}r\right)\d r\, \d z\\
&\lesssim\|r^{-1}\dd_ru_r\|_{L^1(\Omega)}+\|r^{-1}\grad u\|_{L^1(\Omega)}\\
&\lesssim\|r^{-1}\grad u\|_{L^{5/4}(\Omega)}\\
&\lesssim\|\grad u\|_{L^{2}(\Omega)}+\|\grad^2u\|_{L^{5/4}(\Omega)},
}
where we used \eqref{LR_formula} in the third line, and the Hardy inequality \eqref{hardy1} in the last line. 
Next
\eqn{
|I_{2,3,\iin,2}|&=\bigg|\int_{\{r\leq r_0\}} u_\theta \left( \dd_r\frac{u_r}r-\fint_\Omega\p_r \frac{u_r}r\right) \Phi\phi_z'\bigg|\\
&\lesssim\|ru_\theta\|_{L^3(\{r\leq r_0\})}\left\|r^{-\frac12}\left( \dd_r\frac{u_r}r-\fint_\Omega\dd_r\frac{u_r}r\right)\right\|_{L^3(\R^2\times \supp\phi_z)}\|r^{-\frac12}\Phi\|_{L^3(\R^2\times \supp\phi_z)}\\
&\leq A^{O(1)}r_0^\frac23\left\|\grad\dd_r\frac{u_r}r\right\|_{\Luloc^2}\|\grad\Phi\|_{\Luloc^2},
}
where we have used the Hardy inequality (Lemma~\ref{lem_hardy}). Thus Lemma \ref{lem_localgamma} and Young's inequality imply that
\eqn{
\int_{t_0}^t|I_{2,3,\iin,2}|\leq\frac1{20}F(t)+\int_{t_0}^tE.
}
Next let us consider the contributions to $I_2$ from outside $B(r_0)$. Using H\"older's inequality, we obtain that
\eqn{
|I_{2,1,\out}|&=\bigg|\int_{\{r>r_0\}}u_\theta\dd_r\frac{u_r}r\dd_3\Phi\phi_z\,\d x\bigg|\\
&\leq\|u_\theta\|_{\Luloc^6(\{r>r_0\})}\|r^{-1}\dd_ru_r-r^{-2}u_r\|_{\Luloc^3(\{r>r_0\})}\|\grad\Phi\|_{\Luloc^2(\Rd)}.
}
Hence, since Proposition \ref{pointwisebounds} shows that $|u|\leq A^{O(1)} (r^{-1} + r^{-1/4}) $ and $|\p_r u_r | \leq A^{O(1)} (r^{-2} + r^{1/4})$, we see that the first two norms on the right hand side are finite and bounded by, say, $r_0^{-10}$. Thus, an application of Young's inequality gives that
\[\int_{t_0}^t|I_{2,1,\out}|\leq\frac1{20}F(t)+r_0^{-10}(t-t_0).\]
The remaining outer parts of $I_2$, i.e. $I_{2,2,\out}$ and $I_{2,3,\out}$ can be estimated in a similar way, with the latter bounded by, say, $E+r_0^{-10}$.

Finally let us consider $I_3$. Taking $p$ such that, for example, $\frac1p=\frac12-\frac\gamma4$, we have $p-2= 2\gamma /(2-\gamma ) \geq \gamma$, and so our quantified Hardy's inequality (Lemma~\ref{lem_hardy}) shows that
\eqn{
|I_{3,\iin}|&\leq\left\|r^{-2+\frac6p}u_\theta\right\|_{L^{\left(1-\frac2p \right)^{-1}}(\{r\leq r_0\})}\|r^{-\frac3p+\frac12}\Phi\|_{\Luloc^p}\|r^{-\frac3p+\frac12}\Gamma\|_{\Luloc^p}\\
&\lesssim\gamma^{-O(1)}r_0^{\gamma/2}\left( \|\Phi\|_{\Luloc^2}+\|\grad\Phi\|_{\Luloc^2}\right)\left(\|\Gamma\|_{\Luloc^2}+\|\grad\Gamma\|_{\Luloc^2}\right),
}
which gives that $\int_{t_0}^t|I_{3,\iin}|\leq \frac1{20}F(t)+\int_{t_0}^t E$. On the other hand, for $r\geq r_0$ we have the simple bound
\eqn{
|I_{3,\out}|&\leq2 \|r^{-1} u_\theta\|_{L_x^\infty(\{r\geq r_0\})}\|\Phi\|_{\Luloc^2}\|\Gamma\|_{\Luloc^2}\leq r_0^{-5/4}E,
}
as required.\\

\noindent\texttt{Step 3.} Given $\tau >0$ we use the choice of time of regularity (Lemma~\ref{lem_choice}) to find $t_0\in [1-\tau,1]$ such that $E(t_0) \lec A^{O(1)} \tau^{-3}$.\\

Indeed, Lemma~\ref{lem_choice} lets us choose 
 $t_0\in[1-\tau , 1]$ such that 
 \[\|\grad^2 u(t_0)\|_{\infty }\leq A^{O(1)}\tau^{-\frac32}.\]
 It follows from the axial symmetry and \eqref{LR_formula} that $|\Phi | + |\Gamma | \leq | \nabla \omega |$, and so 
\eqnb\label{bd_in_B1}
\| \Phi (t_0) \phi_z^{1/2}  \|_{L^2 (\{ r\leq 1 \} )} + \| \Gamma (t_0) \phi_z^{1/2} \|_{L^2 (\{ r\leq 1 \} )} \lesssim\|\grad\omega(t_0)\|_{L^\infty(B(1)\times\mathbb R)}\leq A^{O(1)}\tau^{-\frac32}
\eqne
for every $z\in \R$. Using the decomposition $\omega=\omega_1^\sharp+\omega_1^\flat$ on the interval $[0,1]$, by \eqref{localenergy}, \eqref{flatbounds}, and H\"older's inequality,
\eqn{
\|\Phi(t_0)\phi_z^{1/2} \|_{L^2( \{ r>1 \} ) }+\|\Gamma(t_0)\phi_z^{1/2}\|_{L^2( \{ r>1 \} ) }&\lesssim\|\omega_1^\sharp\|_{L^2(\Rd)}+\|r^{-1}\omega_1^\flat\|_{L^2 (\{ r>1 \} \cap \supp \, \phi_z ) }\\
&\lesssim \|\grad\unlin_1\|_{L^2(\Rd)}+\|r^{-1}\|_{L_{x'}^4(B(1)^c)}\|\omega_1^\flat\|_{L^4(\Rd)}\\
&\leq A^{O(1)}.
}
This and \eqref{bd_in_B1} proves the claim of this step.\\

\noindent\texttt{Step 4.} We prove the claim.\\

Integration in time of the energy inequality \eqref{en_ineq} from initial time $t_0$ chosen in Step 3 above, taking $\sup_{z\in\mathbb R}$, and applying the estimate \eqref{i1} for $I_1$ and Step 2 for $I_2$, $I_3$ we find that 
\eqn{
E(t)+\frac12F(t)\leq \underbrace{E(t_0)}_{\leq A^{O(1)}\tau^{-3} }+r_0^{-10}+\int_{t_0}^tO(r_0^{-3}+\|u\|_{\infty}+\|\grad^2u\|_{\Lu^{5/4}}+\|\grad u\|_{\Lu^2})E(t')\d t'.
}
for $t\in[t_0,1]$. Thus, by Gr\"onwall's inequality,
\eqn{
E(1)&\leq (A^{O(1)}\tau^{-3}+r_0^{-10})\exp \left( O\left( r_0^{-3}(t-t_0)+A^{O(1)}(t-t_0)^\frac15\right) \right).
}
Setting $\tau\coloneqq r_0^4$, we see that the last exponential function is $O(1)$, and the prefactor gives the required estimate \eqref{claim_of_prop}.
\end{proof}

\section{Proof of Theorem~\ref{regularity}}\label{sec_pf_thm1}

In this section we prove Theorem~\ref{regularity}. Namely, given the $L^{3,\infty}$ bound \eqref{weakL3} on the time interval $[0,1]$, we show that $|\na^j u |\leq \exp \exp A^{O_j (1)}$ at time $1$.\\

\noindent\texttt{Step 1.} We show that $\| b \|_{L^p_{3-\uloc} (\R^3)} \leq C_p \exp \exp A^{O(1)}$ for each $p\in [3,\infty )$, $t\in [1/2,1]$, where $b\coloneqq u_re_r+u_ze_z$ denotes the swirl-free part of the velocity field.\\

To this end we apply Proposition~\ref{prop_energy} to find
\eq{\label{B}
\|\Gamma\|_{L_t^\infty \Luloc^2([\frac12,1]\times\Rd)}\leq \exp\exp A^{O(1)}.
}
On the other hand Proposition \ref{pointwisebounds} shows that
\eqn{
\|r^2\omega\|_{L_x^\infty(\{r\leq10\})}\leq A^{O(1)}.
}
Interpolating between this inequality and \eqref{B} we obtain
\eqn{
\|\omega_\theta\|_{\Luloc^p(\{r\leq10\})}&=\|\Gamma^{\frac23}(r^2\omega_\theta)^{\frac13}\|_{\Luloc^p(\{r\leq10\})}\lesssim\|\Gamma\|_{\Luloc^2}^\frac23\|r^2\omega_\theta\|_{L_x^\infty(\{r\leq10\})}^\frac13\leq\exp\exp A^{O(1)}
}
for all $p\leq3$. Noting that 
\eqn{
\curl b=\omega_\theta e_\theta,\quad\div b=0
}
almost everywhere, and that $\div\, b =0$ we now localize $b$ to obtain an $L^p$ estimate near the axis. Namely, for any unit ball $B\subset\{r\leq10\}$, let $\phi\in C_c^\infty(B)$ such that $\phi\equiv1$ on $B/2$. Observe that for all $p\in[1,3)$ we can use H\"older's inequality for Lorentz spaces \eqref{holder_lorentz} to obtain
\eqn{
\|\div(\phi b)\|_{L^p(\Rd)}=\|b\cdot\grad\phi\|_{p}\lesssim\|b\|_{L^{3,\infty}}\|\grad\phi\|_{L^{3p/(3-p) ,1}}\lesssim A.
}
Applying the Bogovski\u i operator \eqref{bogovskii} to $\div(\phi b)$ on the domain $B\setminus(B/2)$, we find  $\tilde b \in W^{1,p}$ such that $\div\tilde b=0$, $\|b-\tilde b\|_{W^{1,p}(B)}\leq A^{O(1)}$, $\tilde b\equiv b$ in $B/2$, and $\tilde b\equiv0$ outside $B$. Then for any $p\in(1,3)$,
\eqn{
\|b\|_{L^{3p/(3-p)}(B/2)}&\leq\|\tilde b\|_{{3p/(3-p)}}\lesssim\|\grad\tilde b\|_{p}\lesssim\|\hspace{-.025in}\curl\tilde b\|_{L^p(B)}\\
&\leq\|\omega_\theta\|_{L^p(B)}+\|b-\tilde b\|_{W^{1,p}(B)}\\
&\leq \exp\exp A^{O(1)},
}
which is our desired localized estimate. Here we have used the boundedness of the operator $\nabla f \mapsto \curl\, f$ in $L^p$ (which is a consequence of the identity $\curl\,\curl \, f = \grad (\div\, f ) - \Delta f$, which in turn implies that $\nabla f = \nabla  (-\Delta )^{-1} \curl (\curl f) $ for divergence-free $f$). Combining this with the pointwise estimates away from the axis (Proposition~\ref{pointwisebounds}) gives the claim of this step.  \\

\noindent\texttt{Step 2.} We show that there exists $C_0>1$ such that
\eqnb\label{utheta/r12}
\Big\|\frac{u_\theta(t)}{r^{\frac12}}\Big\|_{\Luloc^4}^4\leq\Big\|\frac{u_\theta(t_0)}{r^{\frac12}}\Big\|_{\Luloc^4}^4+1+\exp\exp A^{C_0}\int_{t_0}^t\left\|\frac{u_{\theta }}{r^{\frac12}} \right\|_{L^4_{3-\uloc}}^4
\eqne
for each $t_0 \in [1/2,1]$  and $t\in [t_0,1]$.

To this end we provide a localization of the estimate of $u_\theta/r^{1/2}$ in the spirit of \cite[Lemma~3.1]{cfz_2017}. Indeed, one can calculate from the equation \eqref{utheta_pde} for $u_\theta$ that for a smooth cutoff $\psi=\psi(x_3)$,
\eqn{
&\frac14\frac {\d}{\d t}\int_\Rd\frac{u_\theta^4}{r^2}\psi+\frac34\int_\Rd\Big|\grad\frac{u_\theta^2}r\Big|^2\psi+\frac34\int_\Rd\frac{u_\theta^4}{r^4}\psi_z\\
&\quad=-\frac32\int_\Rd\frac1{r^3}u_ru_\theta^4\psi+\frac18\int_\Rd\frac1{r^2}u_\theta^2(2u_\theta^2u_z-\dd_z(u_\theta^2))\psi'
=:I_1+I_2+I_3.
}

As before, we choose $\psi\in C_c^\infty((-2,2))$ with $\psi\equiv1$ in $[-1,1]$ and define the translates $\psi_z(x):=\psi(x-z)$ for all $z\in\mathbb R$. We consider the energies
\eqn{
E_z(t)\coloneqq \frac14\int_\Rd\frac{u_\theta^4}{r^2}\psi_z ,\qquad & F_z(t)\coloneqq \frac34\int_{t_0}^t\int_\Rd\Big|\grad\frac{u_\theta^2}r\Big|^2\psi_z,\\
E(t)\coloneqq \sup_{z\in\mathbb R}E_z(t),\qquad  &F(t)\coloneqq \sup_{z\in\mathbb R}F_z(t).
}
 By Step 1  and Sobolev embedding,
\eqn{
|I_1|&\lesssim\|u_r\|_{\Luloc^6}\Big\|r^{-\frac12}\frac{u_\theta^2}r\Big\|_{L^{12/5}(\Omega)}^2\\
&\leq \exp\exp A^{O(1)}\left( \Big\|\frac{u_\theta^2}r\Big\|_{L^2(\Omega)}^\frac12\Big\|\grad\frac{u_\theta^2}r\Big\|_{L^2(\Omega)}^\frac32+\Big\|\frac{u_\theta^2}r\Big\|_{L^2(\Omega)}^\frac12\right),
}
where $\Omega \coloneqq \R^2 \times \supp\, \psi$.
It follows that
\eqn{
\int_{t_0}^t|I_1|\leq \frac1{20}F(t)+\exp\exp A^{O(1)}\int_{t_0}^tE+(t-t_0).
}
Similarly,
\eqn{
|I_2|&\lesssim\|u_z\|_{L^6_{3-\uloc}}\Big\|\frac{u_\theta^2}r\Big\|_{L^2_{3-\uloc}}\Big\|\frac{u_\theta^2}r\Big\|_{L^3(\Omega)}\\
&\leq \exp\exp A^{O(1)}E^\frac12\left( \Big\|\frac{u_\theta^2}r\Big\|_{L^2(\Omega)}^\frac12\Big\|\grad\frac{u_\theta^2}r\Big\|_{L^2(\Omega)}^\frac12+\Big\|\frac{u_\theta^2}r\Big\|_{L^2(\Omega)} \right),
}
which yields the same bound as $I_1$. Finally,
\eqn{
|I_3|&=\frac18\bigg|\int_\Rd\frac{u_\theta^2}r\dd_3\frac{u_\theta^2}{r}\psi'\bigg|\lesssim\Big\|\frac{u_\theta^2}r\Big\|_{\Luloc^2}\Big\|\grad\frac{u_\theta^2}r\Big\|_{L^2(\Omega)},
}
so we have
\eqn{
\int_{t_0}^t|I_3|\leq \frac1{20}F(t)+\int_{t_0}^tO(E).
}
Summing and taking the supremum over $z\in\mathbb R$ gives the claim of this step. \\

\noindent\texttt{Step 3.} We deduce that
\eq{\label{done}
\|u\|_{L_t^\infty\Luloc^6([t_0,1]\times\Rd)}\leq \exp\exp A^{O(1)},
}
where 
\[
t_0 \coloneqq 1-\exp(-\exp A^{O(1)}).
\]

Indeed, Lemma~\ref{lem_choice} and Proposition \ref{pointwisebounds} give a $t_0\in[1-\exp(- \exp A^{C_0}),1]$ such that $\|r^{-\frac12}u_\theta(t_0)\|_{L_x^4(\Rd)}\leq \exp\exp A^{2C_0}$. Therefore, applying Gr\"onwall's inequality to the claim of the previous step,
\eqn{
\left\|\frac{u_\theta}{r^{\frac12}}\right\|_{L_t^\infty \Luloc^4([t_0,1]\times\Rd)}\leq\exp\exp A^{O(1)}.
}
Combining this with Proposition \ref{pointwisebounds} and H\"older's inequality,
\eqn{
\|u_\theta\|_{L_t^\infty\Luloc^6([t_0,1]\times\Rd)}&\leq \|ru_\theta\|_{L_x^\infty(\{r\leq1\})}^\frac13\|r^{-\frac12}u_\theta\|_{L_t^\infty \Luloc^4([t_0,1]\times\Rd)}^\frac23+\|u\|_{L_t^\infty L_x^6([t_0,1]\times\{r>1\})}\\
&\leq \exp\exp A^{O(1)},
}
which, together with Step 1, implies \eqref{done}.\\

We note that Step 3 already provides a subcritical local regularity condition of the type of Ladyzhenskaya-Prodi-Serrin, which guarantees local boundedness of all spatial derivatives of $u$, and can be proved by employing the vorticity equation for example (see \cite[Theorem~13.7]{nse_book}).  In the last step below we  use a robust tail estimate of the pressure function (recall Lemma~\ref{lem_poisson_tail}) to provide a  simpler justification of pointwise bounds by $\exp \exp A^{O(1)}$. \\

\noindent\texttt{Step 4.} We prove that, if $\| u \|_{L^\infty ([1-t_1 ,1 ] ; W^{k-1,6}_{\uloc} )} \lec \exp \exp A^{O(1)}$ for some $k\geq 1$ and $t_1=\exp (-\exp A^{O(1)})$, then the same is true for $k$ (with some other $t_1$ of the same order).  \\

Let $I =[a,b] \subset [t_1,1]$, and let $\chi\in C^\infty (\R )$ be such that $\chi (t)=0$ for $t<a+(b-a)/8$ and $\chi(t)=1$ for $t>(a+b)/2$. We set $\phi \in C_c^\infty (B(0,2);[0,1])$ such that $\phi =1$ on $B(0,1/2)$ and $\sum_{j\in \Z^3} \phi_j =1 $, where $\phi_j \coloneqq \phi (\cdot - j)$  for each $j\in \R^3$.

Letting $v\coloneqq \chi \phi \na^k u$ we see that $v(t_1)=0$, and
\[
\begin{split}
v_t -\Delta v &= \underbrace{- \chi' \phi \na^k u - 2\chi \na \phi\cdot \na (\na^k  u) - \chi \Delta \phi (\na^k  u)}_{=: f_1}-\chi \phi \div (1+T )\na^k (u\otimes u ) \\
& = f_1 - \phi \div (1+T ) ((\chi \na^k u \otimes u + u\otimes \chi \na^k u )\tilde \phi )\\
&\qquad -\chi \phi \div(1+T) \sum_{\substack{|\alpha |+|\beta |+|\gamma |=k \\|\alpha |,|\beta |< k}} C_{\alpha,\beta, \gamma}  (D^\alpha u\otimes  D^\beta u D^\gamma \tilde{ \phi } ) -\chi \phi \div T \na^k (u\otimes u (1-\tilde{\phi } ) ) \\
&=: f_1 + f_2 + f_3 + f_4.
\end{split}
\]
We can now estimate $\| v(t) \|_6$, by extracting the same norm on the right-hand side and ensuring that the length of the interval is sufficiently small, so that the norm can be absorbed. Namely, 
\[\begin{split}
\| v(t) \|_{6} &= \left\|\int_{a}^t\ee^{(t-t')\Delta} f_1 (t')\d t' + \int_{a}^t\ee^{(t-t')\Delta}f_2(t' ) \d t'+ \int_{a}^t\ee^{(t-t')\Delta} f_3 (t')\d t' + \int_{a}^t\ee^{(t-t')\Delta} f_4 (t')\d t' \right\|_6 \\
&\leq  \left( \| \chi \na^k u \tilde \phi \|_{L^\infty ([a,t]; L^{6} )} + \| \chi' \na^{k-1} u \tilde\phi \|_{L^\infty ([a,1]; L^{6} )} \right) \int_{a}^t \| \Psi (t-t' ) \|_{W^{1,1}} \d t' \\
&\hspace{2cm}+ \| \chi \na^k u \tilde{ \phi }^{1/2} \|_{L^\infty ([a,t];L^6  )} \| u \tilde{ \phi }^{1/2}\|_{L^\infty ([a,t ]; L^6)  } \int_{a}^t \| \Psi (t-t') \|_{W^{1,6/5}} \d t' \\
&\hspace{2cm}+  \| u \|_{L^\infty ([a,1]; W^{k-1,6}_{\uloc})}^2 \int_{a}^t \| \Psi (t-t' ) \|_{W^{1,6/5}} \d t' \\
&\hspace{2cm}+  \|  \div\, T  (u\otimes u (1-\tilde\phi ) )\|_{L^\infty ([a,1]; W^{k,6}(B(0,2 ))  )} \int_{a}^t \| \Psi (t-t' ) \|_{1} \d t' \\
&\leq \| \chi \na^k u  \|_{L^\infty ([a,t];L^6_{\uloc}  )} \left( (b-a)^{1/2} + \exp \exp A^{O(1)} (b-a)^{1/4} \right) + \exp \exp A^{O(1)}
\end{split}
\]
for each $t\in (a,b)$, where we used Young's inequality, heat estimates \eqref{heat_bounds} and the Calder\'on-Zygmund inequality. By replacing $\phi$ (in the definition of $v$) by $\phi_z$ for any $z\in \R^3$, we obtain the same bound, and so
\[
 \| \chi \na^k u  \|_{L^\infty ([a,b];L^6_{\uloc}  )}\leq  \| \chi \na^k u  \|_{L^\infty ([a,b];L^6_{\uloc}  )} (b-a)^{1/4} \exp \exp A^{O(1)}  + \exp \exp A^{O(1)}.
\]
Thus, for any $b,a$ such that $t_1\leq a<b\leq 1$ and $(b-a)^{1/4} \leq \exp \exp A^{O(1)}/2 $ we can absorb the first term on the right-hand side by the left-hand side to obtain 
\[\| \na^k u \|_{L^\infty ([(a+b)/2,b]; L^6_{\uloc} )} \leq  \exp \exp A^{O(1)}.\]
Since the upper bound is independent of the location of $[a,b]\subset [t_1, 1]$, we obtain the claim.

\section*{Acknowledgements}
WO was partially supported by the Simons Foundation. SP acknowledges support from a UCLA Dissertation Year Fellowship. The authors are grateful to Igor Kukavica, Vladimir \v{S}ver\'ak and Terence Tao for valuable discussions. WO is grateful to Wojciech Zaj{\k a}czkowski for an introduction to the axisymmetric Navier-Stokes equations.

\section*{Conflict of interests statement}
There is no conflict of interests.

\appendix

\section{Quantitative parabolic theory}\label{sec_NU}

Here we prove Proposition~\ref{prop_holder}. Namely, we consider parabolic cylinders
\eqn{
Q_R^{\lambda,\theta}(t_0,x_0)\coloneqq [t_0-\theta R^2,t_0]\times B(x_0,\lambda R),\quad Q_R^{\lambda,\theta}\coloneqq Q_R^{\lambda,\theta}(0,0),\quad Q_R \coloneqq Q_{R}^{1,1}
}
and we consider Lipschitz solutions $V$ of $\mathcal{M}V = 0$ on $Q_R^{\lambda , \theta }$, namely we suppose that
\eqnb\label{MV_eq_weak}
 \int_{\R } \int \left( \p_t V  \phi +\nabla V \cdot \nabla \phi + b\cdot \nabla V \phi \right) =0  
\eqne
for all $\phi \in C_c^\infty (Q_R^{\lambda,\theta} )$, 
where the (distributional) supports of $\div b$ and $V$ are disjoint. Moreover we assume
that \eqref{N} holds, namely
\[
 \mathcal N(R )\coloneqq 2+\sup_{{R'}\leq 2R  }(R')^{-\alpha}\|b\|_{L_t^{\ell}L_x^q(Q_{R'})}<\infty ,
\] 
 where $\alpha\coloneqq \frac nq+\frac2\ell-1\in[0,1)$. We also say that $V$ is a \emph{subsolution} (or \emph{supersolution}) of $\mathcal{M}V=0 $, i.e. $\mathcal{M}V \leq 0$ (or $\mathcal{M}V \geq 0$), if \eqref{MV_eq_weak} holds with ``$=$'' replaced by ``$\leq$'' (or ``$\geq$'') for all nonnegative test functions. 

We will show that 
\eqnb\label{claim_of_prop_holder}
\osc_{B(r)}V(0)\lesssim\left(\frac r R\right)^\gamma\osc_{Q(R)}V
\eqne
for all $r\leq R$, where $\gamma=\exp(-{\mathcal N}^{O(1)})$.\\

To this end we first prove the Harnack inequality for Lipschitz subsolutions of $\mathcal{M} V =0$.

\begin{lem}[based on Lemma 3.1 in \cite{nu_2012}]\label{lem3.1}
Let $V$ be a Lipschitz solution of $\mathcal MV\leq0$ in $Q_R^{\lambda,\theta}$ where $\lambda\in (1,2]$ and $\theta \in (0,1]$.  Then
\eqn{
\sup_{Q_R^{1,\theta/2}}V_+\leq (\mathcal N / \theta )^C\left(\fint_{Q_R^{\lambda,\theta}}V_+^2 \right)^\frac12.
}
\end{lem}

\begin{proof}
We first note that, for any $r,a$ satisfying
\eqn{
\frac 3r+\frac2a \in \left[ \frac 32,\frac 52 \right],
}
 we have the interpolation inequality
\eq{\label{36}
\|\zeta U\|_{L_t^aL_x^r(Q_R^{\lambda,\theta})}&\lesssim_{\lambda,\theta}R^{\frac 3r+\frac2a-\frac 32}\|\zeta U\|_{\mathcal V(Q_R^{\lambda,\theta})},
}
by \cite[(3.4) in Chapter II]{lsu}, where $\mathcal V$ is the energy space $L_t^\infty L_x^2\cap L_{t}^2\dot H^1_x$.\\

Since $V$ is a subsolution, we have, for a non-negative test function $\eta$,
\eqn{
\int_{Q_R^{\lambda,\theta}}(\dd_tV\eta+\grad V\cdot\grad\eta+b\cdot\grad V\eta)\leq0.
}
We let $\eta\coloneqq \varphi'(V)\xi$ where $\xi$ is a cutoff function vanishing on a neighborhood of the boundary of $Q_R^{\lambda,\theta}$, and $\varphi$ is a convex function vanishing on $\mathbb R_-$. Taking $U:=\varphi(V)$ then gives 
\eqn{
\int_{Q_R^{\lambda,\theta}\cap\{V>0\}}\left(\dd_tU\xi+\grad U\cdot\grad\xi+\frac{\varphi''(V)}{\varphi'(V)^2}|\grad U|^2\xi+b\cdot\grad U\xi\right)\leq0.
}
We now take 
\[
\varphi(\tau)\coloneqq \tau_+^p\quad (p>1)\qquad \text{ and }\qquad \xi\coloneqq \chi_{\{t<\overline t\}}U\zeta^2,
\]
where $\zeta$ is a smooth cutoff function in $Q_R^{\lambda,\theta}$ and $\overline t\in(-\theta R^2,0)$,
\eqnb\label{LEI_1}
\int_{B_{\lambda R}}(\zeta U)^2(\overline t)dx+\int_{Q_R^{\lambda,\theta}\cap\{t<\overline t\}}(2-p^{-1})|\grad U|^2\zeta^2+U\grad U\cdot\grad(\zeta^2)+\frac12b\cdot\grad(U^2) \zeta^2-\dd_t(\zeta^2)U^2\leq0.
\eqne
Using integration by parts and recalling the assumption $\div \, b \geq 0$, we can apply H\"older's inequality to obtain 
\eqn{
\int_{Q_R^{\lambda,\theta}\cap\{t<\overline t\}}b\cdot\grad(U^2)\zeta^2&\geq-\int_{Q_R^{\lambda,\theta}\cap-\{t<\overline t\}}b\cdot\grad(\zeta^2)U^2\\
&\hspace{-1in}\geq-\|b\|_{L_t^\ell L_x^q(Q_R^{\lambda,\theta})}\||U|^\frac1s\zeta^{\frac1s-1}\grad\zeta\|_{L_{t,x}^{2s}(Q)}\|(\zeta|U|)^{2-\frac1s}\|_{L_t^{(1-\frac1{2s}-\frac1\ell)^{-1}}L_x^{(1-\frac1{2s}-\frac1q)^{-1}}(Q)}\\
&\hspace{-1in}=-\|b\|_{L_t^\ell L_x^q(Q_R^{\lambda,\theta})}\|U\zeta^{1-s}|\grad\zeta|^s\|_{L_{t,x}^2(Q_R^{\lambda,\theta})}^\frac1s\|\zeta U\|_{L_t^aL_x^r(Q_R^{\lambda,\theta})}^{2-\frac1s}
}
where $s>2$ and $r$ and $a$ are defined by
\eqn{
\frac1{2s}+\frac1q+\frac1r\Big(2-\frac1s\Big)=1,\quad\frac1{2s}+\frac1\ell+\frac1a\Big(2-\frac1s\Big)=1.
}
Applying Young's inequality to separate the last term, and utilizing the interpolation inequality \eqref{36} (which is valid since 
\[
\frac3r+\frac2a=\frac32+1-2\left(1+2/\left(\frac 3q+\frac2\ell\right)\right)^{-1} \in (3/2,11/6),
\]
as needed) we obtain, after plugging into  the local energy inequality \eqref{LEI_1},
\eqn{
\sup_{t\in[-\theta R^2,0]}\int_{B_{\lambda R}}(\zeta U)^2dx+\int_{Q_R^{\lambda,\theta}\cap\{t<\overline t\}}(2-p^{-1})|\grad U|^2\zeta^2+U\grad U\cdot\grad(\zeta^2)-\dd_t(\zeta^2)U^2\\
-O\left( R^2\|b\|_{L_t^\ell L_x^q(Q_R^{\lambda,\theta})}^{2s}\|U\zeta^{1-s}|\grad\zeta|^s\|_{L_{t,x}^2(Q_R^{\lambda,\theta})}^2\right)-\frac1{10}\|\zeta U\|_{\mathcal V(Q_R^{\lambda,\theta})}^2\leq0.
}
Absorbing $\na U$ from the term on the third term on the left-hand side by the second term we obtain 
\eqn{
\|\zeta U\|_{\mathcal V(Q_R^{\lambda,\theta})}^2\lesssim\int_{Q_R^{\lambda,\theta}}\left( |\grad\zeta|^2+\zeta|\dd_t\zeta|+R^2\|b\|_{L_t^\ell L_x^q(Q_R^{\lambda,\theta})}^{2s}\zeta^{2-2s}|\grad\zeta|^{2s}\right) U^2.
}
We now set 
\[ 
\lambda_m\coloneqq 1+2^{-m}(\lambda-1)\qquad \text{ and }\qquad\theta_m\coloneqq \frac12\theta(1+4^{-m}),
\]
and we substitute $\zeta$ with $\zeta_m$ such that
\eqn{
\zeta_m\equiv1\text{  in  }Q_R^{\lambda_{m+1},\theta_{m+1}},\quad\zeta_m\equiv0\text{  outside  }Q_R^{\lambda_m,\theta_m},\quad |\dd_t\zeta_m|\leq \frac{4^mC}{\theta R^2},\quad\frac{|\grad\zeta_m|}{\zeta_m^{1-\frac1s}}\leq \frac{2^mC }{R},
}
where $C$ may depend on $\lambda$. Then the energy estimate and \eqref{36}, taken with $r=l=10/3$, yield
\eqn{
\|\zeta_mU\|_{L_{t,x}^{10/3}(Q_R^{\lambda,\theta})}&\lesssim \|\zeta_mU\|_{\mathcal V(Q_R^{\lambda,\theta})}\leq CR^{-1}(\theta^{-\frac12}+2^m+\mathcal N^s)2^{ms}\|U\|_{L_{t,x}^2(Q_R^{\lambda,\theta})}.
}
Recalling the definition of $U$ and replacing $p$ with $p_m:=(5/3)^m$, H\"older's inequality implies
\eqn{
\left(\fint_{Q_R^{\lambda_{m+1},\theta_{m+1}}}u_+^{2p_{m+1}}\right)^{\frac1{2p_{m+1}}}&\leq\left(C\fint_{Q_R^{\lambda_m,\theta_m}}(\zeta_mU)^{10/3 }\right)^{\frac1{rp_{m}}}\\
&\leq\left(C\theta_m^{-1}\mathcal N^{2s}4^{m(s+1)}\fint_{Q_R^{\lambda_m,\theta_m}}u_+^{2p_m}\right)^{\frac1{2p_m}}.
}
Iterating, we have
\eqn{
\left(\fint_{Q_R^{\lambda_m,\theta_m}}u_+^{2p_m}\right)^{\frac1{2p_m}}&\leq\prod_{k=0}^{m-1}\left( \frac{C}{\theta }4^{k(s+1)}\mathcal N^s\right)^{\frac1{2p_k}}\left(\fint_{Q_R^{\lambda,\theta}}u_+^2\right)^\frac12,
}
and we conclude by taking $m\to\infty$.
\end{proof}

In the next three lemmas we focus on nonnegative solutions to $\mathcal{M}V\leq 0$ and we find lower bounds on the mass distribution of such solutions. We first show that if $V\geq k$ in $Q_R$, except for a small (quantified) ``portion of $Q_R$'', then in fact $V\geq k/2$ everywhere in a smaller cylinder.

\begin{lem}[based on part 2 of Corollary 3.1 in \cite{nu_2012}]\label{cor3.1}
If $V$ is a non-negative solution of $\mathcal MV\geq0$ in $Q_R^{\lambda,\theta}$ and
\eqn{
|\{V<k\}\cap Q_R^{\lambda,\theta}|\leq(\mathcal N/\theta )^{-5C}|Q_R^{\lambda,\theta}|,
}
then
\eqn{
V\geq\frac k2\quad\text{in }Q_R^{1,\theta/2}.
}
\end{lem}

\begin{proof}
We apply Lemma \ref{lem3.1} to $k-V$ to find
\eqn{
\sup_{Q_1^{1,\theta/2}}(k-V)_+\leq(\mathcal N / \theta )^C\left(\fint_{Q_R^{\lambda,\theta}}(k-V)_+^2\right)^\frac12\leq\mathcal N^{-1}k,
}
which implies the result.
\end{proof}

We now show that, if the cylinder $Q_R^{1,\theta }$ is flat enough, then a lower bound on the bottom lid of $Q_R^{1,\theta }$ (i.e. at $t=-\theta R^2$) implies a similar lower bound at every $t$.

\begin{lem}[based on Lemma 3.2 in \cite{nu_2012}]\label{lemma3.2}
Suppose $V$ is non-negative with $\mathcal MV\geq0$ in a neighbuorhood of $Q_R^{1,\theta_0}$ and
\eqn{
|\{V(-\theta_0R^2)\geq k\}\cap B_R|\geq\delta_0|B_R|
}
for some $\delta_0>0$ and $\theta_0\leq C^{-1}\delta_0^6\mathcal N^{-1}$. Then
\eqn{
|\{V(\overline t)\geq\frac13\delta_0k\}\cap B_R|\geq\frac13\delta_0|B_R|
}
for all $\overline t\in[-\theta_0R^2,0]$.
\end{lem}

\begin{proof}
By the calculations  in \cite{nu_2012}, with $\zeta$ a smooth cutoff function supported in $B_R$,
\eq{\label{46}
&\int_{B_R}(V(\overline t)-k)_-^2\zeta^2+\int_{Q_R^{1,\theta_0}}\chi_{\{t<\overline t\}}|\grad(V-k)_-|^2\zeta^2\leq\int_{B_R}(V(-\theta_0R^2)-K)_-^2\zeta^2\\
&\quad\quad+\int_{Q_R^{1,\theta_0}}\chi_{\{t<\overline t\}}(V-k)_-^2\big(O(|\grad\zeta|^2)+b\cdot\grad(\zeta^2)+(\div b)\zeta^2\big).
}
We choose $\zeta$ such that $\zeta\equiv1$ in $B_{(1-\sigma)R}$ and $|\grad\zeta|\leq\frac2{\sigma R}$ where $\sigma<1$ is to be specified. Note that, due to \eqref{div_constraint},
\eqn{
\int_{Q_R^{1,\theta_0}}\chi_{\{t<\overline t\}}(V-k)_-^2(\div b)\zeta^2&\leq k^2\int_{Q_R^{1,\theta_0}}\chi_{\{t<\overline t\}}(\div b)\zeta^2\\
&=-k^2\int_{Q_R^{1,\theta_0}}\chi_{\{t<\overline t\}}b\cdot\grad(\zeta^2).
}
Then the right-hand side of \eqref{46} is bounded by
\eqn{
k^2\Big((1-\delta_0)|B_R|+O(\theta_0\sigma^{-2}|B_R|)+\frac4{\sigma R}\|b\|_{L_t^\ell L_x^q(Q_R)}\|1\|_{L_t^{\ell'}L_x^{q'}(Q_R^{1,\theta_0})}\Big).
}
From here one can proceed with the argument exactly as in \cite{nu_2012} to arrive at
\eqn{
\left| \left\lbrace V(\overline t)<\frac13\delta_0k\right\rbrace \cap B_R\right| \leq \left( 1-\frac13\delta_0\right)^{-2}(1-\delta_0+O(\sigma+\sigma^{-2}\theta_0+\sigma^{-1}\theta_0^{2/\ell'}\mathcal N)).
}
Setting $\sigma=C^{-1/5}\delta_0^2$ and $\theta_0$ as above proves the claimed bound.
\end{proof}
We now show that for any given ``portion of $Q^{1,\theta }_R$'' (in the sense of a set with the measure arbitrarily close to $|Q^{1 , \theta }|$) $V$ is greater or equal a constant multiple of some lower bound, if, for each $t$, the lower bound occurs at least on some ``portion of $B_R$''. Although this enables us to obtain a lower bound on almost the entire cylinder, we lose an exponential in the process.

\begin{lem}[based on Lemma 3.3 in \cite{nu_2012}]\label{lemma3.3}
Let $V\geq 0$ be a solution of $\mathcal MV\geq 0$ in $Q_R^{\lambda,\theta}$ satisfying
\eqn{
|\{V(t)\geq k_0\}\cap B_R|\geq\delta_1|B_R|\text{ for all }t\in[-\theta R^2,0]
}
for some $k_0>0$, $\delta_1>0$. Then for any $\mu>0$ and $s>C(\mathcal N+\theta^{-1})  /(\delta_1\mu)^2$,
\eqn{
|\{V<2^{-s}k_0\}\cap Q_R^{1,\theta}|\leq\mu|Q_R^{1,\theta}|.
}
\end{lem}

\begin{proof}
With $k_m=2^{-m}k_0$, we define
\eqn{
\mathcal E_m(t)\coloneqq \{x\in B_R \colon k_{m+1}\leq V(x,t)<k_m\};\quad\mathcal E_m\coloneqq \{(t,x)\in Q_R^{1,\theta} \colon x\in\mathcal E_m(t)\}.
}
Integrating the inequality $\mathcal MV\geq0$ against the test function $\eta=(V-k_m)_-\xi(x)^2$ where $\xi$ is a smooth cutoff vanishing in a neighborhood of $\partial B_{\lambda R}$ and satisfying $\xi\equiv1$ in $B_R$,
\eq{\label{yes}
\int_{Q^{\lambda,\theta}_{ R} \cap \{ V < k_m \} } | \grad V |^2 \xi^2 &\leq  \int_{Q^{\lambda , \theta }_R } | \grad (V-k_m )_- |^2 \xi^2 \lesssim   \left. \int_{B_{\lambda R}\cap \{ V< k_m \} } (V-k_m)_-^2 \xi^2 \right|_{t=-\theta R^2 } \nonumber\\
&\quad\quad+\int_{-\theta R^2 }^0 \int_{B_{\lambda R}\cap \{ V< k_m \}}(V-k_m)_-^2 |\grad \xi |^2 + 2 (V-k_m)_-^2 \xi b \cdot \grad \xi\nonumber \\
&\lesssim k_m^2R^n(1+\theta\mathcal N),
}
by H\"older's inequality and the trivial bound $0\leq(V-k_m)_-\leq k_m$. From De Giorgi's inequality \cite[(5.6) in Chapter II]{lsu},
\eqn{
(k_m-k_{m+1})|\{V(t)<k_{m+1}\}\cap B_R|\lesssim\frac{R}{\delta_1}\int_{\mathcal E_m(t)}|\grad V(t)|
}
for all $t\in[-\theta R^2,0]$. Integrating in time, squaring, and applying Cauchy-Schwarz gives
\eqn{
k_{m+1}^2\left| \{V<k_{m+1}\}\cap Q_R^{1,\theta}\right|^2\lesssim\frac{R^2}{\delta_1^2}\int_{\mathcal E_m}|\grad V|^2\d x\d t|\mathcal E_m|.
}
Combined with \eqref{yes}, this gives
\eqn{
\left| \{V<k_{m+1}\}\cap Q_R^{1,\theta}\right|^2\lesssim \delta_1^{-2}R^{n+2}(1+\theta\mathcal N)|\mathcal E_m|.
}
We conclude
\eqn{
s\left| \{V<k_s\}\cap Q_R^{1,\theta}\right|^2&\leq\sum_{m=0}^{s-1}\left|\{V<k_{m+1}\}\cap Q_R^{1,\theta}\right|^2\\
&\lesssim\delta_1^{-2}R^{n+2}(1+\theta\mathcal N)\sum_{m=0}^{s-1}|\mathcal E_m|\\
&\lesssim\delta_1^{-2}(\theta^{-1}+\mathcal N)|Q_R^{1,\theta}|^2.
}
\end{proof}

We can now combine Lemmas~\ref{cor3.1}--\ref{lemma3.3} to obtain a pointwise lower bound for $V$ in the interior of a cylinder, with an exponential dependence on $\mathcal{N}$.
\begin{lem}[based on part 1 of Corollary 3.2 in \cite{nu_2012}]\label{cor3.2}
If $V$ is a non-negative solution of $\mathcal MV\geq0$ in $Q_R^{2,1}$ and
    \eqn{
    |\{V(-\Theta R^2)\geq k\}\cap B_R|\geq\delta|B_R|
    }
    for some $k>0$ and $\Theta\leq C^{-1}\delta^6\mathcal N^{-1}$, then
    \eqn{
    V\geq \exp(-\delta^{-2}(\mathcal N /\Theta )^{20C})k\quad\text{in }Q_R^{1,\Theta/2}.
    }
\end{lem}

\begin{proof}
This is a straightforward application of Lemmas \ref{lemma3.2}, \ref{lemma3.3}, and \ref{cor3.1} in sequence, with the latter two applied with $R\to\frac32R$ to compensate for the shrinking domain in Lemma \ref{cor3.1}.
\end{proof}

By considering $V-\inf V$ and $\sup V - V$ the above lemma now allows us to estimate oscillations of solutions to $\mathcal{M} V=0$ with no sign restrictions.

\begin{lem}[based on Lemma 3.5 of \cite{nu_2012}]\label{lemma3.5}
If $V$ solves $\mathcal M V=0$ in $Q_R^{2,1}$ then
\eqn{
\osc_{Q^{(1)}}V\leq(1-\exp(-\mathcal N^{50C}))\osc_{Q^{(2)}}V
}
where $Q^{(1)}=Q_R^{1,\Theta/2}$, $Q^{(2)}=Q_R^{2,1}$, and $\Theta=C^{-2}\mathcal N^{-1}$.
\end{lem}

\begin{proof}
Consider the positive supersolutions $V_1=V-\inf_{Q^{(2)}}V$ and $V_2=\sup_{Q^{(2)}}V-V$. With $k=\osc_{Q^{(2)}}V$, clearly we must have $|\{V_i(-\Theta R^2)\geq k\}\cap B_{2R}|\geq|B_{2R}|/2$ for either $i=1$ or $i=2$. Fix this $i$, so $V_i$ obeys the hypotheses of Lemma \ref{cor3.2}. Let us assume for concreteness that $i=1$; the other case is analogous. Then by the lemma,
\eqn{
\inf_{Q^{(2)}}V+\exp(-\mathcal N^{50C})\osc_{Q^{(2)}}V\leq V\leq\sup_{Q^{(2)}}V
}
for all $(t,x)\in Q^{(1)}$, which immediately implies the result.
\end{proof}

Finally, iterating Lemma~\ref{lemma3.5} we obtain the required H\"older continuity \eqref{claim_of_prop_holder}, i.e. we can prove Proposition~\ref{prop_holder}. 

\begin{proof}[Proof of Proposition~\ref{prop_holder}.]
Iterating Lemma~\ref{lemma3.5}, we have
\eqn{
\osc_{Q^{2,1}_{(\Theta/2)^{k/2}R/2}}V\leq(1-\exp(-\mathcal N^{50C}))^k\osc_{Q_{R/2}^{2,1}}V.
}
We conclude upon taking $k=\lfloor\log\frac Rr(\log\frac2\Theta)^{-1}\rfloor$.
\end{proof}

\bibliographystyle{plain}
\bibliography{literature}

\small
\medskip
\noindent
W.~S.~O\.za\'nski\\
{Department of Mathematics, Florida State University, Tallahassee, FL 32306}\\
e-mail: wozanski@fsu.edu

\medskip
\medskip
\noindent
S.~Palasek\\
{Department of Mathematics, University of California, Los Angeles, CA 90095}\\
e-mail: palasek@math.ucla.edu

\end{document}